\documentclass[11pt]{article}

\usepackage{amsmath, amssymb, amsthm, amsfonts, graphicx} 
\usepackage{cite}
\usepackage{geometry} 

\usepackage{color}

\definecolor{orange}{rgb}{1,0.5,0}

\newtheorem{theorem}{Theorem}[section]

\newtheorem{proposition}[theorem]{Proposition}
\newtheorem{remark}[theorem]{Remark}
\newtheorem{lemma}[theorem]{Lemma}
\newtheorem{corollary}[theorem]{Corollary}

\usepackage{color}

\renewcommand{\AA}{{\mathcal A}}

\newcommand{\EE}{{E}}

\newcommand{\HH}{{\mathcal H}}

\DeclareMathOperator*{\diam}{diam}

\newcommand{\R}{\mathbb{R}}

\newcommand{\eps}{\varepsilon}

\renewcommand{\phi}{\varphi}

\newcommand{\rring}[1]%
{{\mathop{#1}\limits^{\raisebox{-0.2ex}[0ex][0ex]{\tiny{oo}}}}}

\def\XXint#1#2#3{{\setbox0=\hbox{$#1{#2#3}{\int}$}
\vcenter{\hbox{$#2#3$}}\kern-.5\wd0}}

\newcommand{\ncoL}[2]{\|#1\|_{L^\infty({#2})}}

\numberwithin{equation}{section}

\begin{document}

\title{\bf On an isoperimetric problem with a competing non-local
  term. I. The planar case.}

\author{Hans Kn\"upfer\thanks{Hausdorff Center for Mathematics,
    University of Bonn, 53117 Bonn, Germany} \and Cyrill
  B. Muratov\thanks {Department of Mathematical Sciences, New Jersey
    Institute of Technology, Newark, NJ 07102, USA}}
\maketitle
\begin{abstract}
  This paper is concerned with a study of the classical isoperimetric
  problem modified by an addition of a non-local repulsive term. We
  characterize existence, non-existence and radial symmetry of the
  minimizers as a function of mass in the situation where the
  non-local term is generated by a kernel given by an inverse power of
  the distance. We prove that minimizers of this problem exist for
  sufficiently small masses and are given by disks with prescribed
  mass below a certain threshold, when the interfacial term in the
  energy is dominant. At the same time, we prove that minimizers fail
  to exist for sufficiently large masses due to the tendency of the
  low energy configuration to split into smaller pieces when the
  non-local term in the energy is dominant. In the latter regime, we
  also establish linear scaling of energy with mass suggesting that
  for large masses low energy configurations consist of many roughly
  equal size pieces far apart. In the case of slowly decaying kernels
  we give a complete characterization of the minimizers.
\end{abstract}


\section{Introduction}




The isoperimetric problem is a classical problem in the calculus of
variations, one formulation of which seeks to find a set of the
smallest perimeter enclosing a prescribed volume. By the famous result
of De Giorgi, in the Euclidean space $\R^n$ the solution of this
problem is well known to be a ball \cite{degiorgi58}. In this paper we
are interested in the question how the solution of the isoperimetric
problem is affected by an addition of a {\em repulsive long-range
  force}. Specifically, for $n \geq 2$ we wish to study the
variational problem associated with the energy functional
\begin{align} \label{E} %
  E[u] := \int_{\mathbb R^n} |\nabla u| \, dx %
  + \int_{\R^n} \int_{\R^n} {u(x) u(y) \over |x - y|^\alpha} \ dx \,
  dy,
\end{align}
where $u$ is the characteristic function of a subset of $\R^n$ with
finite perimeter and mass $m > 0$, and $\alpha \in (0,n)$ is a
parameter. More precisely, we look for minimizers of the energy $E[u]$
over all $u \in \AA$, where
\begin{align}
  \label{A}
  \AA := \left\{ u \in BV(\R^n, \{ 0, 1 \}) : \int_{\R^n} u \ dx = m
  \right\}.
\end{align}

The choice of the non-local term in (\ref{E}) is motivated by the form
appearing in a number of physical problems. In particular, the
nonlocal term with $\alpha = 1$ in either two or three space
dimensions arises naturally due to Coulombic forces (electrostatic
repulsion in the three-dimensional ambient space), in which case the
characteristic function $u$ may be associated with the uniform charge
density over a subset of either the three-dimensional space or the
two-dimensional plane
\cite{care75,emery93,nagaev95,chen93,nyrkova94,hubert,strukov}. The
case $\alpha = 1$ in three dimensions also arises in the studies of
models of diblock copolymer melts and related polymer, as well as
other, systems (see
e.g. \cite{ohta86,degennes79,stillinger83,glotzer95,m:pre02,ko:book,%
  kovalenko86,mamin94}). More generally, the non-local term is chosen
in this form to have the following four properties \cite{m:phd}:
\begin{itemize}

\item[a)] The non-local term is invariant with respect to translations
  and rotations.

\item[b)] The non-local term is repulsive. 

\item[c)] The non-local term is scale-free.

\item[d)] The non-local term scales with length faster faster than
  volume.
\end{itemize}
Indeed, the kernel in the non-local term depends only on the distance
between two points. The repulsive nature of the non-local term is due
to the fact that for $\alpha > 0$ the kernel is positive and
monotonically decreasing with distance. It is also important to note
that the quadratic form in $L^2(\R^n)$ generated by the kernel is
positive-definite. We point out that when the non-local term has
opposite sign (attractive long-range forces), the minimizers of the
considered variational problem are still balls, since the non-local
term in (\ref{EE}) increases with respect to Schwarz symmetrization
\cite{lieb-loss}. Furthermore, the non-local term is scale-free, in
the sense that dilations only result in the appearance of a
multiplicative factor in front of the non-local term. The scale-free
nature of the non-local term allows one to reduce the number of free
parameters of the problem to a single one, which we choose to be the
mass $m$. Indeed, in this case a coefficient in front of the non-local
term may be eliminated via a rescaling in space (the coefficient in
front of the interfacial term in the energy can be eliminated by
rescaling the energy). Finally, the scaling property is ensured by the
condition $\alpha < n$. Notice that the non-local term is always
infinite when $\alpha \geq n$ and if the interior of the support of
$u$ is non-empty.

We note that a suitably regularized version of the non-local term in
the energy in (\ref{EE}) with $\alpha \in (n, n+1)$ and a negative
sign in front leads to a {\em non-local isoperimetric problem}, in
which the non-local term gives a generalized notion of the perimeter
\cite{caffarelli10}. In fact, when $\alpha \to n + 1$ from below, the
latter $\Gamma$-converges, after a suitable rescaling, to the usual
perimeter \cite{ambrosio11} (see also \cite{caffarelli11}). In our
problem, on the other hand, the non-local term has a very different
effect. It acts as the square of the negative Sobolev norm of $u$ and,
therefore, favors rapid oscillations of $u$. This leads to a {\em
  competition} between the perimeter and the non-local term that can
give rise to the appearance of non-trivial energy minimizing patterns
in bounded domains
\cite{m:phd,m:pre02,m:cmp10,kohn07iciam,choksi10,choksi11,ren07rmp,%
  seul95}. Our whole space problem, in turn, appears as a limit
problem in the studies of $\Gamma$-convergence of the functional in
\eqref{EE} in the presence of a small coefficient multiplying the
perimeter term \cite{choksi10} (see also
\cite{m:cmp10,gms11a,gms11b,choksi11} for a related problem).

Despite the apparent simplicity of the model, for the problem under
consideration even the basic question of existence of minimizers is
not completely straightforward. While in the surface energy-dominated
regime (small masses, $m \ll 1$) one would naturally expect the
minimizers to exist and be in some sense approximations to balls, for
non-local energy-dominated regime (large masses, $m \gg 1$) the energy
may be lowered by splitting a given configuration into several pieces
and moving them far apart. In this situation the minimizers may fail
to exist. Our goal is to address these questions analytically.

In this paper, we present a detailed analysis of existence vs.
non-existence of minimizers of the considered variational problem in
the particular case of two space dimensions. We chose to treat the $n
= 2$ case separately, because in two space dimensions many
technicalities simplify substantially, allowing one to concentrate on
the issues associated with non-locality and making the analysis more
transparent (the general case will be treated elsewhere
\cite{km1}). Furthermore, in two dimensions the obtained results
appear to be optimal, the estimates are readily made explicit, and the
obtained results are applicable to a number of physical systems,
including high-$T_c$ superconductors, magnetic bubble materials and
ferroelectrics \cite{emery93,nagaev95,hubert,strukov}.  What we prove
in the following sections is that for $n = 2$ and $0 < \alpha < 2$ the
basic picture presented above is correct: the minimizer of the
considered problem exists for small enough $m$ and does not exist for
large enough $m$. Note that the considered problem is different from
the one studied in \cite{rigot00}, where the non-local term has a
compactly supported kernel and minimizers exist for all
masses. Moreover, we prove that for $m$ sufficiently small the
minimizer is precisely a single disk.

The main ideas of the proofs are as follows. Existence of minimizers
for small masses is proved by showing that the members of a minimizing
sequence can be chosen to be connected. Non-existence is proved by
showing that for large masses the minimizers must be long and slender,
so it is always possible to reduce the energy by cutting the set in
two and moving the pieces far apart. The fact that the minimizer at
small masses is a ball is proved by exploiting the good stability
properties of the minimizers of the usual isoperimetric problem.

We note that the intricate case of intermediate masses remains largely
open. In particular, an interesting open question is weather the
minimizer of the considered problem is, in fact, a ball whenever it
exists. We prove that this is indeed the case for $\alpha$
sufficiently small. Another interesting open question is about the
structure of the set of values of $m$ for which the minimizer exists,
in particular whether it is an interval. Again, in the case of small
$\alpha$ we prove the latter to be the case and compute the precise
threshold value of $m$ separating the existence and non-existence
regimes. In the full generality, however, these questions are
currently out of reach for the methods of the present paper, since in
our analysis we mainly employ the properties of the problem in the
regimes dominated by either of the two terms in the energy. New tools
that deal with the joint effect of the local and the non-local terms
need to be developed to further address the finer properties of the
minimizers in the considered problem in the regime of intermediate
masses. One step in that direction is the precise scaling of minimal
energy obtained by us for large masses, using an interpolation
inequality relating the interfacial and the non-local parts of the
energy.

Our paper is organized as follows. In Sec. \ref{s:state}, we collect
all the main results of our paper. In Sec. \ref{sec-prelim}, we
present some background results and the results of explicit
computations for several configurations. In
Sec. \ref{sec-exist-minim-small}, we prove existence of minimizers for
sufficiently small masses. In Sec. \ref{sec-non-exist-minim}, we prove
the optimal scaling of the minimal energy for large masses. In
Sec. \ref{sec:non-exist-minim}, we prove non-existence of minimizers
for large masses. In Sec. \ref{sec-shape-minim-small}, we prove that
minimizers for sufficiently small masses are disks. Finally, in
Sec. \ref{sec-alpha}, we prove that for sufficiently small $\alpha$
minimizers exist if and only if they are disks and if and only if
their mass is less or equal than an explicit threshold value.

\section{Statement of results}
\label{s:state}

Throughout the rest of this paper we always assume that $n = 2$ and
$\alpha \in (0, 2)$. The considered variational problem is then
equivalent to minimizing 
\begin{align}
  \label{EE}
  E(\Omega) := |\partial \Omega| + \int_\Omega \int_\Omega {1 \over |x
    - y|^\alpha} \ dx \, dy, \qquad |\Omega| = m,
\end{align}
where $\Omega$ is a set of finite perimeter in the plane, $|\Omega|$
denotes the Lebesgue measure of $\Omega$, i.e., $|\Omega| =
\HH^2(\Omega)$ and $|\partial \Omega|$ denotes the perimeter of
$\Omega$, i.e., $|\partial \Omega| = \HH^1(\partial \Omega)$ (for
definitions, see e.g. \cite{ambrosio}). We note that, when dealing
with minimizers of $E$ in (\ref{EE}), we can always assume that
$\partial \Omega$ is a collection of $C^{2,\beta}$ curves, for some
$\beta \in (0,1)$. This is due to the fact that minimizers of $E$ are
quasi-minimizers of the perimeter and, therefore, the standard
regularity theory of minimal surfaces applies to them. By a
straightforward cutting argument, we also show that minimizers are
connected (but not necessarily simply-connected). The proposition
below collects some basic properties of minimizers:

\begin{proposition} \label{prp-basic} %
  Let $\Omega$ be a minimizer of $E$. Then
  \begin{enumerate}
  \item[(i)] The boundary $\partial \Omega$ of $\Omega$ is of class
    $C^{2,\beta}$, for some $\beta \in (0,1)$, with the regularity
    constants depending only on $m$ and $\alpha$.

  \item[(ii)] $\Omega$ is bounded and connected. Moreover, $\Omega$
    contains at most finitely many holes.

  \item[(iii)] The Euler-Lagrange equation for \eqref{EE} is
      \begin{align} \label{EL}
        \kappa(x) + 2 v(x) - \mu = 0, \qquad v(x) := \int_\Omega {1 \over |x
          - y|^\alpha} \ dy,
      \end{align}
     where $\kappa(x)$ and $v(x)$ are the curvature and the non-local potential
      at the point $x \in \partial \Omega$, respectively, and $\mu \in \R$ is
      the Lagrange multiplier due to the mass constraint (the sign of the
      curvature is chosen to be positive for convex sets).
  \end{enumerate}
\end{proposition}

\noindent Note that further $C^{3,\beta}$ regularity can be inferred
when $\alpha < 1$, since $v \in C^{1,\beta}(\R^2)$ for some $\beta \in
(0,1)$ in that case.

We now state our main results. Concerning the existence of minimizers,
as was noted already in the introduction the question is not
straightforward, since the minimizing sequences for (\ref{EE}) may
consist of disconnected pieces moving off to infinity away from each
other. What we can establish, however, is that in the regime of small
masses, i.e., when the perimeter is the dominant term in the energy,
minimizers of $E$ do indeed exist.

\begin{theorem}[Existence of minimizers] \label{thm-exist} %
 There is $m_1 = m_1(\alpha) > 0$ such that for all $m \leq m_1$
  there exists a minimizer of $E$.
\end{theorem}

Our next result gives a complete characterization of the minimizers
for sufficiently small masses.

\begin{theorem}[Disk as a unique minimizer] \label{thm-disk} %
  There is $m_0 = m_0(\alpha) > 0$ such that for all $m \leq m_0$ the
  unique, up to translations, minimizer of $E$ is $\Omega = B_R(0)$
  with $R = (m / \pi)^{1/2}$.
\end{theorem}

\noindent Since in the case of small masses the energy is dominated by
the perimeter, it is expected that the minimizer should be close to a
disk. However, our result is stronger, stating that even for small
(but positive) masses the minimizers is precisely a disk with mass
$m$. At the same time, our next theorem shows that the minimizer
(global or local) cannot be a disk if $m$ is sufficiently large. Note
that these kinds of results have been derived for a number of related
problems
\cite{mo1:pre96,m:pre02,choksi11,goldstein96,seul95,hubert,thiele70}.

\begin{theorem}[Global and local instability of a
  disk] \label{thm-inst} %
  Let $\Omega = B_R(0)$ with $R = (m / \pi)^{1/2}$. There are two
  constants $m_{c1} < m_{c2}$ given in \eqref{mc1} and \eqref{mc2}
  such that the following holds:
  \begin{itemize}
  \item[(i)] $\Omega$ is not a global minimizer if $m > m_{c1} =
    m_{c1}(\alpha) > 0$.
  \item[(ii)] $\Omega$ is not a local minimizer (with respect to
    arbitrarily small perturbations of the boundary) if $m > m_{c2} =
    m_{c2}(\alpha) > 0$.
  \end{itemize}
\end{theorem}

\noindent Note that by Theorem \ref{thm-inst} disks cease to be global
minimizers before they undergo shape instability. While Theorem
\ref{thm-inst} only shows that a minimizer cannot be given by a disk
with mass $m$ if $m$ is sufficiently large, the next theorem
encompasses a more general result: For sufficiently large masses the
energy \eqref{EE} does not have a minimizer. This point was
conjectured in \cite[Remark 4.2]{choksi10} for the case $n = 3$ and
$\alpha = 1$. In fact, a more general non-existence result holds for
all $\alpha < 2$ in dimensions $n > 2$ as well \cite{km1}.

\begin{theorem}[Non-existence] \label{thm-nonexistence} %
  There exists $m_2 = m_2(\alpha) > 0$ such that there is no
  minimizer of $E$ for all $m > m_2$.
\end{theorem}

\noindent The non-attainability of the minimal energy for large masses
established by Theorem \ref{thm-nonexistence} is related to the fact
that for $m \gg 1$ it is advantageous for the mass to escape to
infinity. One could imagine that in this case the minimizing sequence
consists asymptotically of disconnected sets of approximately equal
masses $m_i \sim 1$ (for a related result in two dimensions, see
\cite{m:cmp10,choksi10,gms11a}). In particular, the minimal energy
should scale linearly with the mass $m$ for large $m$. The next
theorem supports this picture.
\begin{theorem}[Scaling and equipartition of
  energy] \label{thm-linear} %
  Let $m_0 = m_0(\alpha) > 0$ be as in Theorem \ref{thm-disk}. Then
  there exist two constants $C,c > 0$ only depending on $\alpha$ such
  that for all $m > m_0$ we have
  \begin{align} \label{EcC}%
    c m \ \leq \ \inf E(\Omega) \ \leq \ C m.
  \end{align}
  Furthermore, we have equipartition of energy in the sense that for
  every $m \geq m_0$ and every configuration $\Omega$ with $E(\Omega)
  \ \leq C m$, both terms in the energy obey the same bounds as in
  \eqref{EcC} separately.
\end{theorem}

\noindent Note that Theorem \ref{thm-linear} shows in particular that
for large $m$ the minimal scaling of the energy can only be reached if
the interfacial and the nonlocal part of the energy are of the same
order. The result in Theorem \ref{thm-linear} is a consequence of the
multiplicative interpolation inequality \eqref{int}, which is derived
in the course of the proof of Lemma \ref{lin-lindown}.

Theorems \ref{thm-exist}, \ref{thm-disk} and \ref{thm-nonexistence}
cover the two extremes of the range of values of $m$, but say nothing
about what happens at intermediate values of $m$. Thus, the global
structure of minimizers for all masses is currently not
available. Nevertheless, when $\alpha$ is sufficiently small, i.e.,
when the non-local interaction is slowly decaying with distance, we
have a complete characterization of minimizers:

\begin{theorem}[Complete characterization of minimizers for slowly
  decaying kernels] 
  \label{thm-al0}
  Let $m_{c1} = m_{c1}(\alpha) > 0$ be given by \eqref{mc1}.  There
  exists a universal constant $\alpha_0 > 0$ such that for all $\alpha
  \leq \alpha_0$ we have:
  \begin{enumerate}
  \item[(i)] For all $m \leq m_{c1}$ there exists a minimizer of $E$;
    this minimizer, up to translations, is given by $B_R(0)$ with $R =
    (m / \pi)^{1/2}$.
  \item[(ii)] For all $m > m_{c1}$ there is no minimizer of $E$.
  \end{enumerate}
\end{theorem}

\noindent In other words, for sufficiently small values of $\alpha$
the constants in Theorems \ref{thm-exist}--\ref{thm-nonexistence} obey
$m_0 = m_1 = m_2 = m_{c1}$. Thus, our results support a recent
conjecture of \cite{choksi11} in the present setting for sufficienly
slowly decaying kernels. We note that the arguments in the proof of
Theorem \ref{thm-al0} also imply that as soon as it is known that for
a given value of $\alpha$ the minimizers of $E$ can only be disks, we
have $m_0 = m_1 = m_2 = m_{c1}$ without the need to assume that
$\alpha$ is small. The latter follows from a global result in Lemma
\ref{lem-inst-1}(ii) concerning minimizers of $E$ among sets
consisting of two arbitrary disks.



\section{Preliminaries}
\label{sec-prelim}

In this section, we present the necessary background for the proofs of
the main theorems, as well as the results of some precise ansatz-based
calculations. We first outline the proof of the regularity of
$\partial \Omega$ for minimizers:

\begin{lemma}
  \label{lin-qm}
  Let $\Omega$ be a minimizer of $E$.  Then $\partial \Omega$ is of
  class $C^{2,\beta}$, for some $\beta \in (0,1)$, with the regularity
  constants depending only on $m$ and $\alpha$.
\end{lemma}

\begin{proof}
 Regularity of $\partial \Omega$ follows from the fact that every minimizer of
  $\EE$ is a quasiminizer of the interfacial energy. More precisely, we claim
  that for any set of finite perimeter $\Omega' \subset \R^2$ with $|\Omega'| =
  m$, we have
  \begin{align}
    \label{qm}
    |\partial \Omega| \leq |\partial \Omega'| + c (2 -
    \alpha)^{-{\alpha \over 2}} m^{2 - \alpha \over 2} |\Omega \Delta
    \Omega'|,
  \end{align}
  for some universal $c > 0$\footnote{Here and throughout the rest of
    the paper $A \Delta B := (A \backslash B) \cup (B \backslash A)$
    denotes the symmetric difference of the sets $A$ and
    $B$.}. Indeed, assuming that \eqref{qm} holds, we can immediately
  apply \cite[Theorem 1.4.9]{rigot00} to conclude the uniform
  $C^{1,\beta}$ regularity of $\partial \Omega$ for any minimizer
  $\Omega$ of $E$. Furthermore, the boundary satisfies the weak form
  of the Euler-Lagrange equation \eqref{EL}. Noting that clearly $v
  \in C^{0,\beta}(\R^2)$ for some $\beta \in (0,1)$ and since
  $\partial\Omega$ is locally a graph of a $C^{1,\beta}$ function, the
  regularity assertion of the lemma follows by further application of
  the regularity theory for graphs \cite{giusti} (for details of the
  argument, see e.g. the last paragraph of Sec. 2 in
  \cite{sternberg11}, as well as \cite{gilbarg}).

  It hence remains to show \eqref{qm}. By a direct computation, for
  any $R > 0$ we have
  \begin{align}
    \label{dEnloc}
    \hspace{6ex} & \hspace{-6ex} \left| \int_{\Omega'} \int_{\Omega'}
      {1 \over |x - y|^\alpha} \ dx \, dy - \int_\Omega \int_\Omega {1
        \over |x - y|^\alpha} \ dx \, dy \right| \ %
    \leq \ 2 \int_{\Omega \Delta \Omega'} \int_{\Omega \cup \Omega'}
    {1 \over |x - y|^\alpha} \ dx
    \, dy  \nonumber \\
    &\leq \ 2 |\Omega \Delta \Omega'| \sup_{x \in \R^2} \left(
      \int_{B_R(x)} {1 \over |x - y|^\alpha} \ \, dy + \int_{(\Omega
        \cup \Omega') \backslash B_R(x)} {1 \over |x - y|^\alpha} \ \,
      dy \right) \nonumber \\
    &\leq \ 2 |\Omega \Delta \Omega'| \left( {2 \pi R^{2-\alpha} \over
        2 - \alpha} + m R^{-\alpha} \right).
  \end{align}
  The statement then follows from the minimizing property of $\Omega$
  by choosing $R = (2 - \alpha)^{1/2} m^{1/2}$ in (\ref{dEnloc}).
\end{proof}

\begin{lemma}
  \label{prp-connect}
  Let $\Omega$ be a minimizer of $E$. Then $\Omega$ is connected.
\end{lemma}

\begin{proof}
  By uniform H\"older estimates of Lemma \ref{lin-qm}, we have $\Omega
  \subset B_{R_0}(0)$ for some $R_0 > 0$ (where $R_0$ depends on the
  configuration). Therefore, if $\Omega^{(1)}$ is a connected
  component and $\Omega \backslash \Omega^{(1)} \not= \varnothing$,
  then for any $R > 0$ we can consider a set $\Omega'$ obtained by
  translating $\Omega^{(1)}$ outside $B_{R_0 + R}(0)$. The energy of
  the obtained set is then
  \begin{align}
    \label{EomiR}
    E(\Omega') & \leq E(\Omega) - 2 \int_{\Omega^{(1)}} \int_{\Omega
      \backslash \Omega^{(1)}} {1 \over |x - y|^\alpha} \ dx \, dy + 2
    R^{-\alpha} |\Omega^{(1)}| (m - |\Omega^{(1)}|) \notag \\
    & \leq E(\Omega) + 2 (R^{-\alpha} - (2 R_0)^{-\alpha})
    |\Omega^{(1)}| (m - |\Omega^{(1)}|).
  \end{align}
 But for $R$ sufficiently large, this inequality contradicts the minimizing
  property of $\Omega$.
\end{proof}

\begin{lemma} \label{lem-topo} %
  Let $\Omega$ be a minimizer of $E$. Then $\Omega$ contains at most
  finitely many holes.
\end{lemma}

\begin{proof}
 By Lemma \ref{lin-qm} we can cover $\partial \Omega$ by finitely
  many balls of equal radius, which depends only on $\alpha$ and $m$,
  such that in each ball $\partial \Omega$ is a graph of a $C^1$
  function. Therefore, since the perimeter of $\Omega$ is bounded,
  $\partial \Omega$ breaks into a finite collection of simple closed
  curves.
\end{proof}

We now turn to the exact computations related to sets enclosed by
ellipses. In the lemma below, we obtain an expression for the energy
of such a set.

\begin{lemma}
 \label{lin-ellipse}
  Let $\Omega_e$ be a set enclosed by an ellipse of eccentricity
  $e$. Then
  \begin{align}
    \label{Eellipse}
    E(\Omega_e) & = {4 R \over \sqrt[4]{1 - e^2}} \mathrm{E}(e^2) +
    \frac{\pi^2 (1 - e^2)^{-\frac{\alpha+2}{4}} \Gamma (2-\alpha) R^{4
        - \alpha}}{ \Gamma \left(2-\frac{\alpha}{2}\right) \Gamma
      \left(3-\frac{\alpha}{2}\right)} \nonumber \\
    & \times \left\{ (1 - e^2) \,
      _2F_1\left(\frac{1}{2},1-\frac{\alpha}{2};1;e^2 \right) +
      (1-e^2)^{\alpha/2} \,
      _2F_1\left(\frac{1}{2},1-\frac{\alpha}{2};1;\frac{e^2}{e^2 -
          1}\right)\right\},
  \end{align}
  where $\mathrm{E}(x)$ is the complete elliptic integral of the
  second kind, $_2F_1(a, b; c; z)$ is the hypergeometric function, and
  $R = (m / \pi)^{1/2}$.
\end{lemma}

\begin{proof}
  We consider the region $\Omega_e$ enclosed by an ellipse whose
  semi-axes are $a$ and $b$. Since the area of the ellipse is $\pi a b
  = \pi R^2$, we have $a = R / \sqrt[4]{1 - e^2}$ and $b = R
  \sqrt[4]{1 - e^2}$. We also recall that the perimeter of the ellipse
  is given by the well-known expression
  \begin{align}
    \label{EOmeper}
    |\partial \Omega_e| \ = \ {4 R \over \sqrt[4]{1 - e^2}}
    \mathrm{E}(e^2),
  \end{align}
  where $\mathrm{E}(x)$ is the complete elliptic integral of the
  second kind \cite{abramowitz}. To compute the non-local part
  $E_{nl}$ of the energy, we pass to the Fourier space. In terms of
  the Fourier transform
  \begin{align}
   \label{uq}
   \hat u_q = \int_{\R^2} e^{i q \cdot x} u(x) \, dx,
  \end{align}
  of $u = \chi_{\Omega_e}$, the characteristic function of $\Omega_e$,
  the nonlocal part of the energy is given by (see
  e.g. \cite{lieb-loss})
  \begin{align} \label{EOmenonq} %
    E_{nl}(\Omega_e) := \int_{\R^2} \int_{\R^2} {u(x) u(y) \over |x -
      y|^\alpha} \, dx dy = {1 \over (2 \pi)^2} \int_{\R^2} |\hat
    u_q|^2 G_q \, d q,
  \end{align}
  where 
  \begin{align}
    \label{Gq}
    G_q = {2^{2-\alpha} \pi \Gamma \left( 1 - {\alpha \over 2} \right)
      \over \Gamma \left( {\alpha \over 2} \right)} \, |q|^{\alpha-2},
  \end{align}
  is the Fourier transform of the kernel in the non-local term and
  $\Gamma(x)$ is the Gamma-function \cite{abramowitz}. To proceed, we
  note that if $x = (x_1, x_2)$, the rescaling $x_1 \to x_1/a$ and
  $x_2 \to x_2 / b$ transforms $u(x)$, after a suitable translation,
  into the characteristic function of $B_R(0)$. Therefore, with $q =
  (q_1, q_2)$ one can write explicitly upon integration
 \begin{align}
   \label{uq12}
   \hat u_q = 2 \pi R \left( q_1^2 \sqrt{1-e^2} + {q_2^2 \over \sqrt{1
         - e^2}} \right)^{-1/2} J_1 \left( R \sqrt{ q_1^2 \sqrt{1-e^2}
       + {q_2^2 \over \sqrt{1 - e^2}} } \right) ,
 \end{align}
 where $J_1(x)$ is the Bessel function of the first kind
 \cite{abramowitz}. Performing another change of variables: $\tilde q_1
 = q_1 \sqrt[4]{1 - e^2}$ and $\tilde q_2
 = q_2 / \sqrt[4]{1 - e^2}$, and then introducing polar coordinates
 $\tilde q_1 = s \cos t$, $\tilde q_2 = s \sin t$, upon integration in
 $s$ we obtain
 \begin{align}
   \label{EOmenlt}
   E_{nl}(\Omega_e) = {2^{2 - \alpha \over 2} \pi (1 - e^2)^{2 -
       \alpha \over 4} \Gamma(2 - \alpha) R^{4 - \alpha} \over \Gamma
     \left( 2 - {\alpha \over 2} \right) \Gamma \left( 3 - {\alpha
         \over 2} \right)} \int_0^{2 \pi} (2 - e^2 + e^2 \cos 2
   t)^{\alpha - 2 \over 2} dt.
 \end{align}
 Finally, performing the integration in $t$, we obtain
 \begin{align}
   \label{EOmeEnl}
   E_{nl}(\Omega_e) & = \frac{\pi^2 (1 - e^2)^{-\frac{\alpha+2}{4}}
     \Gamma (2-\alpha) R^{4 - \alpha}}{ \Gamma
     \left(2-\frac{\alpha}{2}\right) \Gamma
     \left(3-\frac{\alpha}{2}\right)} \nonumber \\
   & \times \left\{ (1 - e^2) \,
     _2F_1\left(\frac{1}{2},1-\frac{\alpha}{2};1;e^2 \right) +
     (1-e^2)^{\alpha/2} \,
     _2F_1\left(\frac{1}{2},1-\frac{\alpha}{2};1;\frac{e^2}{e^2 -
         1}\right)\right\},
 \end{align}
 where $_2F_1(a, b; c; z)$ is the hypergeometric function
 \cite{abramowitz}. Combining this formula with (\ref{EOmeper}), we
 obtain the result.
\end{proof}

Setting $e = 0$, we also obtain the precise formula for the energy of
a single ball of mass $m = \pi R^2$.
\begin{corollary}
  \label{cor-ball}
  We have 
  \begin{align}
    \label{EBR}
    E(B_R(0)) = 2 \pi R + {2 \pi^2 \Gamma(2 - \alpha) \over \Gamma
      \left( 2 - {\alpha \over 2} \right) \Gamma \left( 3 - {\alpha
          \over 2} \right) } \, R^{4 - \alpha},
  \end{align}
\end{corollary}

The result in Lemma \ref{lin-ellipse} enables us to give the proof of
Theorem \ref{thm-inst} concerning the failure of minimality (either
global or local) of $\Omega = B_R(0)$ for $R = (m / \pi)^{1/2}$ large
enough. The proof relies on the explicit formula for the energy of
elliptical domains (balls, in particular) obtained in Lemma
\ref{lin-ellipse} and Corollary \ref{cor-ball}. The proof of Theorem
\ref{thm-inst} follows directly from the next two lemmas.  Note that
it is easy to see that $m_{c1} < m_{c2}$ where $m_{c1}$ and $m_{c2}$
are defined in \eqref{mc1} and \eqref{mc2}.
\begin{lemma} \label{lem-inst-1} %
  Let
  \begin{align}
    \label{mc1}
    m_{c1}(\alpha) = \pi \left( \frac{\left(\sqrt{2}-1\right) \Gamma
        \left(2-\frac{\alpha}{2}\right) \Gamma
        \left(3-\frac{\alpha}{2}\right)}{\pi \left(1 - 2^{\alpha - 2
            \over 2}\right) \Gamma (2-\alpha)} \right)^{2 \over 3 -
      \alpha}.
  \end{align}
  Then
  \begin{enumerate}
  \item[(i)] The disk of area $m$ is not a global minimizer of $E$ if
    $m > m_{c1}$.
  \item[(ii)] The disk of area $m$ has lower energy than any two
    non-overlapping disks of the same total area if $m \leq m_{c1}$.
  \end{enumerate}
\end{lemma}

\begin{proof}
  We first compare the energy of $\Omega = B_R(0)$ with that of
  $\Omega' = B_{R/\sqrt{2}} (r e_1) \cup B_{R / \sqrt{2}} (-r e_1)$,
  where $R = (m / \pi)^{1/2}$, $r > 2 R$, and $e_1$ is the unit vector
  along the $x_1$-axis. It is easy to see from \eqref{EBR} that
  \begin{align}
    \label{EEpdiskglob}
    E(\Omega') - E(\Omega) \leq 2 \pi (\sqrt{2} - 1) \, R + {\pi^2
      \Gamma(2 - \alpha) \over \Gamma \left( 2 - {\alpha \over 2}
      \right) \Gamma \left( 3 - {\alpha \over 2} \right) } \,
    (2^{\alpha/2} - 2) \, R^{4 - \alpha} + {2^{-1-\alpha} m^2 \over (r
      - R)^\alpha} < 0,
  \end{align}
  if $m > m_{c1}$, where $m_{c1}$ is defined in \eqref{mc1} and $r$ is
  sufficiently large, contradicting minimality of $\Omega$.

  We now show that it is energetically advantageous to replace a set
  $\Omega$ consisting of two non-overlapping disks of mass $tm$ and
  $(1 - t)m$, with arbitrary $t \in (0, 1)$, by a single disk of mass
  $m$ whenever $m \leq m_{c1}$. Indeed, by positivity of the kernel in
  the non-local part of the energy and \eqref{EBR} we have $ E(\Omega)
  > 2 \sqrt{\pi} \, m^{1/2} \tilde F(t, m/m_{c1})$, where
  \begin{align}
    \label{EBt}
    \tilde F(t, \mu) := t^{1/2} + (1 - t)^{1/2} + {2 (\sqrt{2}-1)
      \over 2 - 2^{\alpha/2}} \mu^{3 - \alpha \over 2} \left( t^{4 -
        \alpha \over 2} + (1 - t)^{4 - \alpha \over 2} \right), \qquad
    \mu := {m \over m_{c1}}.
  \end{align}
  The statement of the lemma is equivalent to showing that $F(t, \mu)
  := \tilde F(t, \mu) - \tilde F(0, \mu) \geq 0$ for all $\mu \leq 1$
  and, by symmetry, for all $t \in [0,1/2]$.  We claim that for fixed
  $t$ the minimum of $F$ as a function of $\mu$ is attained when $\mu
  = 1$. Indeed, differentiating this expression with respect to $\mu$,
  we find that
  \begin{align}
    \label{mu}
    {\partial F \over \partial \mu} = {\sqrt{2}-1 \over 2 -
      2^{\alpha/2}} (3 - \alpha) \mu^{1 - \alpha \over 2} F_1(t),
    \qquad F_1(t) := t^{4 - \alpha \over 2} + (1 - t)^{4 - \alpha
      \over 2} - 1.
  \end{align}
  Clearly, $F_1(t)$ is strictly convex and in view of the fact that
  $F_1(0) = 0$, we have $F_1(t) < 0$ and, hence, $\partial F
  / \partial \mu < 0$ for all $t \in (0,1/2]$. Therefore, it is
  sufficient to show that $F(t) \geq 0$ for $\mu = 1$.

  We now write
  \begin{align}
    \label{FF}
    F(t, 1) = {2(\sqrt{2}-1) \over 2 - 2^{\alpha/2}} F_1(t) + F_2(t),
    \qquad F_2(t) := t^{1/2} + (1 - t)^{1/2} - 1.
  \end{align}
  To prove that this function is non-negative for all $\alpha \in
  (0,2)$ and $t \in [0,1/2]$ is a tedious exercise in calculus (the
  reader can readily verify this fact by plotting $F$ as a function of
  two variables, $\alpha$ and $t$). Below we sketch the analytical
  argument. Introduce a cutoff parameter $t_0 = 1/8$. It is then not
  difficult to see that for all $t \in [t_0, 1/2]$ we have
  \begin{align}
    F_1(t) \geq & \ 2^{\alpha - 2 \over 2} - 1 + 2^{\alpha - 6 \over
      2} \left(2^\alpha 7^{-\alpha/2}+1\right)
    \left(1-\frac{\alpha}{2}\right) \left(2-\frac{\alpha}{2}\right) (1
    - 2 t)^2, \\
    F_2(t) \geq & \ \sqrt{2} - 1 - {4 \sqrt{2} \over 25} (1 - 2 t)^2.
  \end{align}
  Then, by inspection we find that $F(t, 1) \geq 0$ for all $t \in
  [t_0, 1/2]$ and all $\alpha \in (0,2)$.

  We now turn to $t \in (0, t_0)$. Here we have the following
  estimates:
  \begin{align}
   F_1(t) \ \geq \ 2^{\alpha - 2 \over 2} - 1 - \left( 2 - {\alpha
        \over 2} \right) \, t, &&
    F_2(t) \ \geq \ \sqrt{2} - 1 + {\sqrt{2} \, ( 2 \sqrt{7} - 1)
      \over \sqrt{7}} \, t.
  \end{align}
  Again, by inspection these estimates imply that $F(t, 1) \geq 0$ for
  all $t \in (0, t_0)$ and all $\alpha \leq 1$. To cover the range
  $\alpha > 1$, we use a different estimate for $F_1$:
  \begin{align}
    F_1(t) \geq \frac12 ( 2 - \alpha ) \big( t \ln t + (1 - t) \ln ( 1
    - t) \big).
  \end{align}
  After some further manipulations, one can obtain a lower bound
  \begin{align}
    F(t, 1) \geq \left(2 \sqrt{2} -\frac{\left(\sqrt{2}-1\right) (2-
        \alpha) (1+\ln 7)}{2 - 2^{\alpha/2}} -\sqrt{\frac{2}{7}} \,
    \right) t,
  \end{align}
  which implies that $F(t, 1) \geq 0$ for all $t \in (0, t_0)$, when
  $\alpha > 1$.
\end{proof}

\begin{lemma} \label{lem-inst-2} %
  Let
  \begin{align}
    \label{mc2}
    m_{c2}(\alpha) = \pi \left( \frac{3 \Gamma
        \left(2-\frac{\alpha}{2}\right) \Gamma
        \left(3-\frac{\alpha}{2}\right)}{\pi \alpha \Gamma (3-\alpha)}
    \right)^{2 \over 3 - \alpha}.
  \end{align}
 Then the disk of area $m$ is not a local minimizer of $E$ (with respect to
  arbitrarily small perturbations of the boundary) if $m > m_{c2}$.
\end{lemma}

\begin{proof}
  We expand the expression in \eqref{Eellipse} for the energy of
  $\Omega_e$ as in Lemma \ref{lin-ellipse} in the power series in $e$
  at $e = 0$:
  \begin{align}
    \label{DEl2} 
    E(\Omega_e) - E(B_R(0)) = \left( {3 \pi R \over 32} - {\pi^2
        \alpha \Gamma(3 - \alpha) R^{4 - \alpha} \over 32 \Gamma
        \left( 2 - {\alpha \over 2} \right) \Gamma \left( 3 - {\alpha
            \over 2} \right) } \right) e^4 + O(e^6) < 0,
  \end{align}
  for sufficiently small $e$ if $m > m_{c2}$, where $m_{c2}$ is
  defined in \eqref{mc2}. Under this condition the energy decreases
  upon arbitrarily small distortion of a disk into an ellipse.
\end{proof}

Lastly, we also need to characterize the non-local potential generated
by a unit ball and, specifically, its behavior near the boundary, in
order to prove Theorem \ref{thm-disk}.

\begin{lemma}
 \label{lin-vball}
  Let 
  \begin{align}
    \label{vB}
    v^B(x) := \int_{B_1(0)} {1 \over |x - y|^\alpha} \, dy.
  \end{align}
  Then
  \begin{align}
    \label{vbr}
    v^B(x) = 
    \begin{cases}
      \left( {\pi \over |x|^\alpha} \right) \,
      _2F_1\left(\frac{\alpha}{2}, \frac{\alpha}{2};
        2;\frac{1}{|x|^2}\right), & |x| \geq 1, \vspace{1mm} \\
      \vspace{1mm} \left( {2 \pi \over 2-\alpha} \right) \,
      _2F_1\left(\frac{\alpha-2}{2}, \frac{\alpha}{2}; 1;
        |x|^2\right), & |x| < 1.
    \end{cases}
  \end{align}
  where $_2F_1(a, b; c; z)$ is the hypergeometric function. In
  particular, if $r = |x| - 1$, we have
  \begin{align}
    \label{vbrser}
    v^B(x) - v_0 =
    \begin{cases}
      - \frac{\pi \alpha (2-\alpha) \Gamma (1-\alpha)}{2 \Gamma^2
        \left(2-\frac{\alpha}{2}\right)} \, r +O(|r|^{2 - \alpha}), &
      \alpha < 1, \\
      -r (2 \ln |r|^{-1} -2 + 3 \ln 4)+O(r^2 \ln |r|^{-1}), & \alpha
      = 1, \\
      -\frac{\sqrt{\pi} \, \Gamma \left(\frac{\alpha-1}{2}\right)
      }{(2-\alpha) \Gamma \left(\frac{\alpha}{2} \right)} \,
      |r|^{1-\alpha} r + O(r), & \alpha > 1,
    \end{cases}
  \end{align}
  where $\Gamma(x)$ is the Gamma-function and
  \begin{align}
    \label{v0}
    v_0 := \frac{\pi \Gamma (2-a)}{\Gamma^2
      \big(2-\frac{a}{2}\big)}.
  \end{align}
\end{lemma}

\begin{proof}
  The proof is by an explicit computation. Introducing the Fourier
  transform $\hat v_q^B$ of $v^B$:
  \begin{align}
    \label{vbq}
    \hat v^B_q = \int_{\R^2} e^{i q \cdot x}
    v^B(x) \, dx,
  \end{align}
  and using (\ref{Gq}) and (\ref{uq12}) with $e = 0$, we obtain
  \begin{align}
    \label{vbqq}
    \hat v^B_q = G_q \hat u_q = {2^{3-\alpha} \pi^2 \Gamma \left( 1 -
        {\alpha \over 2} \right) \over \Gamma \left( {\alpha \over 2}
      \right)} \, |q|^{\alpha-3} J_1(|q|),
  \end{align}
  where $\hat u_q$ is the Fourier transform of the characteristic
  function of the unit ball centered at the origin.  

  Inverting the Fourier transform and integrating over the directions
  of $q$, with $z = |q|$ we arrive at
  \begin{align}
    \label{vbqqq}
    v^B(x) = {1 \over (2 \pi)^2} \int_{\R^2} e^{-i q \cdot x} \hat
    v^B_q \, dq = {2^{2 - \alpha} \pi \Gamma \left( 1 - {\alpha \over
          2} \right) \over \Gamma \left( {\alpha \over 2} \right) }
    \int_0^\infty z^{\alpha - 2} J_1(z) J_0(z |x|) dz,
  \end{align}
  where $J_n(x)$ are the Bessel functions of the first kind. But the
  right-hand side of (\ref{vbqqq}) coincides with the right-hand side
  of (\ref{vbr}). Finally, the expansion in (\ref{vbrser}) is an
  immediate consequence of (\ref{vbr}).
\end{proof}

\section{Existence of minimizers for small masses}
\label{sec-exist-minim-small}

We now prove the existence result in Theorem \ref{thm-exist}. The
strategy of the proof is to suitably localize the minimizing sequence
for $E$ in \eqref{EE}. Existence of minimizers then follows by the
usual compactness and lower-semicontinuity results for functions of
bounded variation \cite{ambrosio}.

\begin{proof}[Proof of Theorem \ref{thm-exist}]
  Let $\{ \Omega_k \}_{k=1}^{\infty}$, with $\Omega_k \subset \R^2$
  and $|\Omega_k| = m$, be a minimizing sequence for $E$. Without loss
  of generality, we may assume that each $\Omega_k$ consists of $N_k <
  \infty$ disjoint open connected components $\Omega_k^{(i)}$ ordered
  so that $|\Omega_k^{(1)}| \geq |\Omega_k^{(2)}| \geq \ldots \geq
  |\Omega_k^{(N_k)}|$, and that the interfaces $\partial\Omega_k$ are
  smooth. As a first step, we use a ball of radius $R = (m /
  \pi)^{1/2}$ as a test function to obtain an upper bound for the
  minimal energy. By comparing the energy of $\Omega$ with the energy
  of the ball $B_R(0)$, we may assume that
  \begin{align} \label{EOmn3al} %
    E(\Omega_k) \ \leq \ E(B_R(0)) \ = \ %
    2 \sqrt{\pi} \, m^{1/2} (1 + C m^{3 - \alpha \over 2}),
  \end{align}
  for some $C > 0$ depending only on $\alpha$ (for the precise
  constant, see \eqref{EBR}).

  Suppose now that $N_k > 1$ and, hence, $|\Omega_k^{(i)}| \leq m/2$
  for all $2 \leq i \leq N_k$. By the isoperimetric inequality and by
  positivity of the non-local term in the energy, we have
  \begin{align}
    \label{EisoBR}
    2 \sqrt{\pi} \, (|\Omega_k^{(i)}|^{1/2} + (m -
    |\Omega_k^{(i)}|)^{1/2}) \leq E(\Omega_k) \qquad \text{for all } 2
    \leq i \leq N_k.
  \end{align}
  Squaring both sides of (\ref{EisoBR}) and combining it with
  (\ref{EOmn3al}), after some algebraic manipulations we obtain 
  \begin{align}
    \label{Omcm4al}
    |\Omega_k^{(i)}| \leq C m^{4-\alpha} \qquad \text{if} \qquad 2
    \leq i \leq N_k \quad \text{and} \quad m \leq 1,
  \end{align}
  for some $C > 0$ depending only on $\alpha$.

  On the other hand, consider a set $\Omega_k'$ obtained by erasing
  $\Omega_k^{(N_k)}$ from $\Omega_k$ and then dilating the resulting
  set by $\lambda_k = \sqrt{ m / (m - |\Omega_k^{(N_k)}|)} \in (1,
  \sqrt{2}]$, so that $|\Omega_k'| = m$ once again. If $E(\Omega_k') <
  E(\Omega_k)$, we replace the set $\Omega_k$ by $\Omega_k'$ in the
  minimizing sequence and repeat the above process. Then, after
  finitely many steps either the set $\Omega_k$ is connected, or
  $E(\Omega_k') \geq E(\Omega_k)$. In the latter case we can write
  \begin{align} \label{Eomdil} %
    E(\Omega_k') \ &= \ \lambda_k (|\partial \Omega_k| - |\partial
    \Omega_k^{(N_k)}|) + \lambda_n^{4-\alpha} \int_{\Omega_k
      \backslash \Omega_k^{((N_k)}} \int_{\Omega_k \backslash
      \Omega_k^{(N_k)}} {1 \over |x - y|^\alpha} \ dx \,
    dy \nonumber \\
    &\leq \lambda_k^4 E(\Omega_k) - |\partial \Omega_k^{(N_k)}|.
  \end{align}
  and, therefore,
  \begin{align}
    \label{dOmtoosmall}
   |\partial \Omega_k^{(N_k)}| \leq \ (\lambda_k^4 - 1) E(\Omega_k) \ %
    \leq \ \frac{6 |\Omega_k^{(N_k)}|}m \ E(\Omega_k).
  \end{align}
  Applying again the isoperimetric inequality on the left-hand side of
  \eqref{dOmtoosmall} and using the fact that by \eqref{EOmn3al} we
  have $E(\Omega_k) \leq Cm^{1/2}$ for some $C > 0$ depending only on
  $\alpha$ and $m \leq 1$, we then conclude that in this case
  \begin{align}
    \label{Omitoosmall}
    |\Omega_k^{(i)}| \geq c m \qquad \text{if} \qquad 1 \leq i \leq
    N_k \quad \text{and} \quad m \leq 1,
  \end{align}
  for some $c > 0$ depending only on $\alpha$. 

  It is easy to see that for sufficiently small $m$ the two
  inequalities in (\ref{Omcm4al}) and (\ref{Omitoosmall}) are
  incompatible. Thus, given a minimizing sequence
  $\{\Omega_k\}_{k=1}^\infty$, for sufficiently small $m$ it is always
  possible to construct another minimizing sequence
  $\{\Omega_k'\}_{k=1}^\infty$, in which each set $\Omega_k'$ is
  connected.

  By a suitable translation, one can further assume that the origin
  belongs to each of $\Omega_k'$. In turn, since the perimeters of
  $\Omega_k'$ are uniformly bounded above, we have $\Omega_k' \Subset
  B_R(0)$ for some large enough $R > 0$. Therefore, introducing the
  characteristic functions $u_k \in BV(B_R(0); \{0, 1\})$ of
  $\Omega_k'$, we get that the functions $u_k$ are equibounded in
  $BV(B_R(0); \{0, 1\})$. So up to extraction of a subsequence $u_k
  \to u \in BV(B_R(0); \{0, 1\})$ strongly in $L^1(B_R(0))$ and $u_k
  \rightharpoonup u$ in $BV(B_R(0); \{0, 1\})$, with the limit
  independent of $R$. In particular, $\int_{B_R(0)} u \, dx =
  m$. Since the perimeter is lower-semicontinuous, and the non-local
  term is continuous with respect to the above convergence (the latter
  follows immediately from (\ref{dEnloc})), we conclude that the set
  $\Omega = \{u = 1\}$ is a minimizer.
\end{proof}

\section{Scaling of the minimal energy}
\label{sec-non-exist-minim}

We now consider the opposite extreme, in which the non-local term
favors splitting of the set $\Omega$ into smaller disconnected
sets. The corresponding scaling of the minimal energy is described by
Theorem \ref{thm-linear} whose proof is an immediate consequence of
the following three lemmas. We note that the main point of Theorem
\ref{thm-linear} is the ansatz-free lower bound for large $m$ which
matches the upper bound from an ansatz consisting of a collection of
equal size balls far apart. We also note that we only need to prove
the bounds in Theorem \ref{thm-linear} for sufficiently large
masses. Indeed, by the isoperimetric inequality and by the positivity
of the non-local term we have $E(\Omega) \geq |\partial \Omega| \geq 2
\sqrt{\pi} \, m^{1/2}$, so $E(\Omega)$ is uniformly bounded away from
zero whenever $m \geq c$, for any $c > 0$.

We begin with an ansatz-based upper bound.

\begin{lemma}
 \label{lin-linup}
 For every $m \geq 1$ there exists $\Omega$ such that $E(\Omega) \leq
 C m$ for some $C > 0$ depending only on $\alpha$.
\end{lemma}

\begin{proof}
  The proof is by an explicit construction. We take
  \begin{align}
    \label{OmBRN}
    \Omega \ = \ \left(\bigcup_{n = 0}^{N-1} B_1(n R e_1) \right) \cup
    B_r(NR e_1),
  \end{align}
  where $N = \lfloor (m / \pi) \rfloor$, $e_1$ is the unit vector
  along the $x_1$ direction, $r = \pi^{-1/2} (m - \pi N)^{1/2}$, and
  $R > 2$, i.e., we take $\Omega$ to be a linear chain of
  non-overlapping unit balls (except for the last one, whose radius is
  chosen to accommodate the mass constraint). Then by \eqref{EBR} we
  have
  \begin{align} \label{EBmN} %
    E(\Omega) \ \leq \ (N+1)E(B_1(0)) + {2 \pi^2 N(N+1) \over (R -
      2)^\alpha} \ \leq \ C m + {4 m^2 \over (R -
      2)^\alpha}.
  \end{align}
  So the assertion of the lemma follows by choosing $R = m^{1/\alpha}
  + 2$.
\end{proof}

We now turn to the ansatz-free lower bound.

\begin{lemma}
 \label{lin-lindown}
 For every admissible $\Omega$ we have $E(\Omega) \geq c m$, for some
 universal $c > 0$.
\end{lemma}

\begin{proof}
  The proof of the lower bound can be obtained by retracing the steps
  in the proof of \cite[Lemma B.1]{m:cmp10}. Here we present a simpler
  proof, which does not rely on Fourier techniques and the properties
  of special functions. The result is obtained from the following
  interpolation inequality:
  \begin{align} \label{int} %
    \int_{\R^2} u^2 \, dx \ \leq \ C \left( ||u||_{L^\infty(\R^2)}
      \int_{\R^2} |\nabla u| \ dx \right)^{\frac{2-\alpha}{3-\alpha}}
    \left( \int_{\R^2} \int_{\R^2} \frac{u(x)u(y)}{|x-y|^\alpha} \ dx
      dy \right)^{\frac{1}{3-\alpha}},
  \end{align}
  for some universal $C > 0$, which is valid for any $u \in BV(\R^2)
  \cap L^\infty(\R^2)$. Indeed, for any admissible set $\Omega$, let
  $u$ be the characteristic function of $\Omega$. Applying
  \eqref{int}, we then have 
  \begin{align}
    \label{lbest7}
    m & \leq C \left(\int_{\R^2} |\nabla u| \ dx
    \right)^{\frac{2-\alpha}{3-\alpha}} \left( \int_{\R^2} \int_{\R^2}
      \frac{u(x)u(y)}{|x-y|^\alpha} \ dx dy
    \right)^{\frac{1}{3-\alpha}} \notag \\
    & \leq C E^{\frac{2-\alpha}{3-\alpha}}(\Omega)
    E^{\frac{1}{3-\alpha}}(\Omega) = C E(\Omega).
  \end{align}

  For the proof of \eqref{int}, we only need to take into account the
  non-local interaction on intermediate length scales of order $R >
  0$, which will be determined later. With the change of variables $z
  = y-x$, we have
  \begin{align}
    \hspace{6ex} & \hspace{-6ex} \int_{\R^2} \int_{\R^2}
    \frac{u(x)u(y)}{|x-y|^\alpha} \ dx dy \ \geq \
    \int_{\R^2}\int_{B_{2R}(0) \backslash B_R(0)}
    \frac{u(x)u(x+z)}{|z|^\alpha} \ dz dx \notag \\
    &= \ \int_{\R^2} \int_{B_{2R}(0) \backslash B_R(0)}
    \frac{|u(x)|^2}{|z|^\alpha} \ dz dx %
    + \int_{\R^2} \int_{B_{2R}(0) \backslash B_R(0)} \frac{u(x)
      (u(x+z)-u(x))}{|z|^\alpha} \ dz dx.
  \end{align}
  Using the fact that $R \leq |z| \leq 2R$ and that $|B_{2R}(0)
  \backslash B_R(0)| = 3\pi R^2$, we hence get
  \begin{align}
    \hspace{2ex} & \hspace{-2ex}
    \int_{\R^2} \int_{\R^2} \frac{u(x)u(y)}{|x-y|^\alpha} \ dx dy \notag \\
    &\geq \ C R^{2-\alpha} \int_{\R^2} u^2 \, dx - R^{-\alpha}
    \int_{\R^2} \int_{B_{2R}(0) \backslash B_R(0)} \int_0^1
    |z| |u(x)| |\nabla u(x + tz)| \ dt dz dx \notag \\
    &\geq \ C R^{2-\alpha} \int_{\R^2} u^2 \ dx - C' R^{3-\alpha}
    ||u||_{L^\infty(\R^2)} \int_{\R^2} |\nabla u(x)| \ dx,
  \end{align}
  for some universal $C, C' > 0$ (recall that $\alpha \in (0,2)$),
  where we argued by approximating $u$ with smooth functions, noting
  that by an argument similar to the one used in \eqref{dEnloc} the
  non-local term is continuous in the $L^1$-topology. Therefore
  \begin{align} \label{oben} %
    \int_{\R^2} u^2 \, dx \ & \leq \ C R ||u||_{L^\infty(\R^2)}
    \int_{\R^2} |\nabla u| \ dx + C R^{\alpha-2} \int_{\R^2}
    \int_{\R^2} \frac{u(x)u(y)}{|x-y|^\alpha} \ dx dy,
  \end{align}
  for some universal $C > 0$.  Estimate \eqref{int} then follows by
  minimizing the right hand side of \eqref{oben} in $R$, i.e. we
  choose
  \begin{align}
    R \ = \ \left( ||u||_{L^\infty(\R^2)} \int_{\R^2} |\nabla u| \ dx
    \right)^{-\frac 1{3-\alpha}} \left( \int_{\R^2} \int_{\R^2}
      \frac{u(x)u(y)}{|x-y|^\alpha} \ dx dy \right)^{\frac
      1{3-\alpha}}.
  \end{align}
  This concludes the proof of \eqref{int}.
\end{proof}

The following lemma strengthens the lower bound for configurations
which satisfy the linear scaling of the energy:
\begin{lemma} \label{lin-equi} Let $m \geq 1$ and suppose that
  $\Omega$ satisfies $E(\Omega) \leq C m$ for some $C > 0$.  Then
  there is another constant $c > 0$ depending only on $C$ such that
  \begin{align}
    |\partial \Omega| \geq c m, && \text{and} && \int_{\Omega}
    \int_{\Omega} \frac 1{|x-y|^\alpha} \, dx dy \geq c m.
  \end{align}
\end{lemma}
\begin{proof}
  Follows directly from $E \leq Cm$ and \eqref{int}.
\end{proof}

\section{Non-existence of minimizers for large masses}
\label{sec:non-exist-minim}

We now present the proof of Theorem \ref{thm-nonexistence}. We begin
with a basic estimate (from above and below) of the diameter of a
minimizer $\Omega$ of $E$.

\begin{lemma}
 \label{lin-d}
 Let $m \geq 1$, let $\Omega$ be a minimizer of $E$ and let $d :=
 \diam (\Omega)$. Then
  \begin{align}
    \label{dm}
    c m^{1/\alpha} \leq d \leq C m,
  \end{align}
  for some $c, C > 0$ depending only on $\alpha$.
\end{lemma}

\begin{proof}
  We first recall that by Proposition \ref{prp-basic} the set $\Omega$
  is regular and connected. Therefore, we have $2 d \leq |\partial
  \Omega| \leq E(\Omega) \leq C m$ for some $C > 0$ depending only on
  $\alpha$, in view of Lemma \ref{lin-linup}.

  On the other hand,
  \begin{align}
    \label{Enllb}
    {m^2 \over d^\alpha} \leq \int_\Omega \int_\Omega {1 \over |x -
      y|^\alpha} \ dx \, dy \leq E(\Omega) \leq C m,
  \end{align}
  which yields the second inequality.
\end{proof}

Let us note that as an immediate consequence of Lemma \ref{lin-d} we get
non-existence of minimizers for large masses when $\alpha < 1$. We
call this regime {\em far field-dominated}, as opposed to the opposite
regime ($\alpha \geq 1$), which we call {\em near-field dominated}. We
will see again in Sec. \ref{sec-shape-minim-small} that this
distinction also plays a role for minimizers at small masses.

\begin{corollary}[Non-existence in the far field-dominated case]
  \label{c:nonlong}
  Let $\alpha < 1$. Then there exists $m_2 = m_2(\alpha) > 0$ such
  that the there are no minimizers of $E$ for all $m > m_2$.
\end{corollary}

\begin{proof}
  By Lemma \ref{lin-linup} and \ref{lin-d} every minimizer $\Omega$
  has to satisfy $E \leq C m$ and $E \geq \diam(\Omega) \geq c
  m^{1/\alpha}$. For $\alpha < 1$ and sufficiently large $m$, both
  inequalities cannot be satisfied at the same time.
\end{proof}

We now turn to completing the proof of Theorem \ref{thm-nonexistence},
which in view of Corollary \ref{c:nonlong} amounts to the proof of the
following proposition.

\begin{proposition}[Non-existence in the near field-dominated case]
 \label{prp-nonnear}
 Let $\alpha \geq 1$. Then there exists $m_2 = m_2(\alpha) > 0$ such
 that there are no minimizers of $E$ for all $m > m_2$.
\end{proposition}

\begin{proof}
  We argue by contradiction. Let $\Omega$ be a minimizer of $E$ for
  some $m \geq 1$. Introducing $d := \diam(\Omega)$, let $x_1, x_2
  \in \partial \Omega$ be such that $|x_1 - x_2| = d$. For every $s
  \in (0,d)$, define $T(s)$ to be the line perpendicular to $x_1 -
  x_2$ and located at distance $s$ from $x_1$. The line $T(s)$ cuts
  the set $\Omega$ into two non-empty parts. We define the part of
  $\Omega$ which is closer to $x_1$ as $\Omega_s$. We also define
  $V(s) := |\Omega_s|$ and $A(s) := |\Omega \cap T(s)|$. Note that $A
  \in L^\infty(0,d)$ (since the diameter of $\Omega$ is bounded), and
  by Cavalieri's principle we have
  \begin{align}
    \label{VAs}
    V(s) = \int_0^s A(s') \, ds' \qquad \forall s \in (0, d).
  \end{align}
  Also, without loss of generality we may assume that $V(d/2) \leq
 m/2$. In particular, this implies that for all $s \in (0,d/2)$,
  \begin{align}
    \label{EOmnon}
    E(\Omega) \geq |\partial \Omega| + \int_{\Omega_s} \int_{\Omega_s}
    {1 \over |x - y|^\alpha} \ dx \, dy + \int_{\Omega \backslash
      \Omega_s} \int_{\Omega \backslash \Omega_s} {1 \over |x -
      y|^\alpha} \ dx \, dy + {m V(s) \over d^\alpha}.
  \end{align}
  
  Now, consider a new set $\Omega' = (\text{T}_R \Omega_s) \cup
  (\Omega \backslash \Omega_s)$, where $\text{T}_R$ denotes a
  translation by distance $R > 0$ along the vector
  $\overrightarrow{x_2 x_1}$, i.e., the set $\Omega'$ is obtained by
  cutting $\Omega$ with $T(s)$ and moving the resulting pieces
  distance $R$ apart. We have
  \begin{align}
    \label{EOmpnon}
    E(\Omega') & \leq |\partial \Omega| + 2 A(s) \notag \\
    & + \int_{\Omega_s} \int_{\Omega_s} {1 \over |x - y|^\alpha} \ dx
    \, dy + \int_{\Omega \backslash \Omega_s} \int_{\Omega \backslash
      \Omega_s} {1 \over |x - y|^\alpha} \ dx \, dy + {2 m V(s) \over
      R^\alpha}.
  \end{align}
  Therefore, from the minimizing property of $\Omega$ we obtain
  \begin{align}
    \label{dEnon}
    2 A(s) \geq {m V(s) \over 2 d^\alpha} \qquad \forall s \in (0,
    d/2),
  \end{align}
  for large enough $R$, or, equivalently,
  \begin{align}
    \label{dEnon2}
    {d V \over ds} \geq { m V \over 4 d^\alpha} \qquad \text{for
      a.e. } s \in (0, d/2).
  \end{align}
  Integrating this expression from $s \in (0, d/2)$ to $d/2$, we then
  conclude that
  \begin{align}
    \label{Vexp}
    V(s) \leq {m \over 2} \exp \left( -{ m (d - 2 s) \over 8 d^\alpha}
    \right) \qquad \forall s \in (0, d/2).
  \end{align}
 In particular, by Lemma \ref{lin-d} 
  \begin{align}
    \label{Vd4}
    V(s) \leq {m \over 2} \exp \left( - \frac{1}{16} m d^{1 - \alpha}
    \right) \leq m \exp \left( -c m^{2 - \alpha} \right) \qquad
    \forall s \in (0, d/4],
  \end{align}
  for some $c > 0$ depending only on $\alpha$, i.e., $V(s)$ becomes
  uniformly small for $s \in (0, d/4]$ and $m \gg 1$.

  Let us now show that the latter is impossible. We consider a
  different set $\Omega''$ obtained by erasing $\Omega_s$ from
  $\Omega$ and dilating the resulting set $\Omega \backslash \Omega_s$
  by a factor $\lambda_s = \sqrt{m / (m - V(s))} > 1$ to make
  $\Omega''$ admissible. By the minimizing property of $\Omega$ and
  positivity of the kernel in the non-local term, we have
  \begin{align}
    \label{EOmppnon}
    E(\Omega) \leq E(\Omega'') \leq \lambda_s^{4-\alpha} E(\Omega
    \backslash \Omega_s) \leq \lambda_s^4 (E(\Omega) - |\partial
    \Omega_s| + 2 A(s)),
  \end{align}
  where we argued as in \eqref{Eomdil}. Therefore, by isoperimetric
  inequality and Lemma \ref{lin-linup} we have
  \begin{align}
    \label{dEOmppnon}
    2 \sqrt{\pi} \, V^{1/2}(s) \leq |\partial \Omega_s| \leq C (V(s) +
    A(s)),
  \end{align}
  for some $C > 0$ depending only on $\alpha$. In view of (\ref{Vd4}),
  there exists $m_2 \geq 1$ such that $C V \leq \sqrt{\pi} \, V^{1/2}$
  for all $s \in (0, d/4]$ and all $m > m_2$. Therefore, for these
  values of $m$ (\ref{dEOmppnon}) implies
  \begin{align}
    \label{dVpp}
    {dV \over ds} \geq c V^{1/2} \qquad \text{for a.e. } s \in
    (0, d/4),
  \end{align}
  with some $c > 0$ depending only on $\alpha$. Integrating this
  inequality from 0 to $s \in (0, d/4]$, we then find that
  \begin{align}
    \label{Vsmnon}
    V(s) \geq c s^2 \qquad \forall s \in (0, d/4],
  \end{align}
 for some $c > 0$ depending only on $\alpha$.  But by Lemma \ref{lin-d}
  this contradicts (\ref{Vd4}) at $s = d/4$.
\end{proof}

\section{Shape of minimizers for small masses}
\label{sec-shape-minim-small}

We now turn to the proof of Theorem \ref{thm-disk}. Here it is
convenient first to rescale length in such a way that the rescaled set
$\Omega$ has a fixed mass. Let us define a positive parameter
\begin{align} \label{def-eps} %
  \eps := \left( {m \over \pi} \right)^{3-\alpha \over 2}.
\end{align}
Then the renormalized energy
\begin{align}
  \label{Eeps}
  E_\eps(\Omega) := |\partial \Omega| + \eps \int_\Omega \int_\Omega
  {1 \over |x - y|^\alpha} \ dx \, dy, \qquad |\Omega| = \pi,
\end{align}
is related to the original energy as 
\begin{align}
  \label{Erenorm}
  E(\Omega) = (m / \pi)^{1/2} E_\eps(\Omega_\eps),
\end{align}
where $\Omega_\eps$ is obtained by dilating the set $\Omega$ by a
factor of $(m / \pi)^{-1/2}$. We note that by virtue of Theorem
\ref{thm-exist}, the minimizers of $E_\eps$ exists for all $\eps \leq
\eps_1(\alpha)$, where $\eps_1$ is related to $m_1$ via
(\ref{def-eps}). Furthermore, the regularity result in Proposition
\ref{prp-basic}, with constants depending on $\eps$ and $\alpha$,
holds for the minimizers of $E_\eps$.

Expressed in terms of the rescaled problem, Theorem \ref{thm-disk}
takes the following form:
 \begin{proposition}
   There exists $\eps_0 = \eps_0(\alpha) > 0$ such that for all $\eps
   \leq \eps_0$ the minimizer of $E_\eps$ is a unit disk.
 \end{proposition}

 \noindent The proof proceeds differently for the far field-dominated
 ($\alpha < 1$) and near field-dominated ($\alpha \geq 1$)
 regimes. But before we turn to the proof, let us establish a number
 of basic properties of the minimizers of $E_\eps$ that we will need
 in our analysis. Recall that for any set of finite perimeter in
 $\R^2$, the isoperimetric deficit is given by
 \begin{align} \label{isodef} %
   D(\Omega) := {|\partial \Omega| \over 2 \pi} - 1.
 \end{align}
 We begin with a basic estimate of the isoperimetric deficit of
 minimizers.

 \begin{lemma} \label{lin-isodef} %
   Let $\Omega$ be a minimizer of $E_\eps$, and let $D(\Omega)$ be the
   isoperimetric deficit of $\Omega$. Then for some $C > 0$ depending
   only on $\alpha$ we have
   \begin{align}
     \label{Deps0}
    D(\Omega) \leq C \eps.
   \end{align}
\end{lemma}

 \begin{proof}
  The proof is obtained by testing $E_\eps$ with a unit disk. The assertion
   follows immediately from the minimizing property of $\Omega$, positivity of
   the non-local term and the fact that by \eqref{EBR},
   \begin{align}
     \label{dOmepsD}
    |\partial \Omega| \leq E_\eps(\Omega) \leq E_\eps(B_1(0)) = 2 \pi + C \eps.
   \end{align}
\end{proof}

 We next establish that for small values of $\eps$ the minimizers are
 necessarily convex and, hence, simply connected.

 \begin{lemma}
  \label{lin-convex}
   Let $\Omega$ be a minimizer of $E_\eps$. Then there exists $\eps_2
  = \eps_2(\alpha) > 0$ such that $\Omega$ is convex for all $\eps
   \leq \eps_2$.
 \end{lemma}

 \begin{proof}
   The Euler-Lagrange equation for the minimizers of $E_\eps$ is given
   by
   \begin{align}
     \label{ELeps}
     \kappa(x) + 2 \eps v(x) - \mu = 0, \qquad v(x) := \int_\Omega {1
       \over |x - y|^\alpha} \ dy,
   \end{align}
   where, as in \eqref{EL}, $\kappa(x)$ and $v(x)$ denote the
   curvature and potential at $x \in \partial \Omega$, respectively,
   and $\mu \in \R$ is the Lagrange multiplier. To estimate $\mu$, let
   us integrate (\ref{ELeps}) over the outer boundary $\partial
   \Omega_o$ of $\Omega$, which is justified by Proposition
   \ref{prp-basic}. After dividing by $|\partial \Omega_o| > 0$, we
   obtain
   \begin{align}
     \label{mueps}
     \mu = {2 \pi \over |\partial \Omega_o|} + {2 \eps \bar v}, \qquad
     \bar v := {1 \over |\partial \Omega_o|} \int_{\partial \Omega_o}
     v(x) \, d \HH^1(x).
   \end{align}
   Now, let $\Omega_o$ be the set enclosed by $\partial \Omega_o$, so
   that $\Omega \subseteq \Omega_o$ and $\partial \Omega_o
   \subseteq \partial \Omega$. In particular, we have $\pi = |\Omega|
   \leq |\Omega_o|$, and by the isoperimetric inequality $2\sqrt{\pi}
   \, |\Omega_o|^{1/2} \leq |\partial \Omega_o| \leq |\partial
   \Omega|$.  Lemma \ref{lin-isodef} thus implies
   \begin{align}
     \label{dOmoo} 2\pi \leq |\partial \Omega_o| \leq 2 \pi + C \eps,
   \end{align} 
   for some $C > 0$ depending only on $\alpha$.  Similarly, we have
   for every $x \in \Omega$
   \begin{align}
     \label{vest}
     0 \leq v(x) \leq \int_{B_1(x)} {1 \over |x - y|^\alpha} \, dy \,
     + \int_{\Omega \backslash B_1(x)} {1 \over |x - y|^\alpha} \, dy
    \leq C && \text{and} && 0 \leq \bar v \leq C,
   \end{align}
   for some $C > 0$ depending only on $\alpha$. Inserting
   \eqref{dOmoo} and \eqref{vest} into \eqref{mueps}, we obtain that
   $|\mu - 1| \leq C \eps$ for some $C > 0$. Substituting this
   estimate, together with (\ref{vest}), into (\ref{ELeps}), we then
   conclude that $|\kappa(x) - 1| \leq C \eps$. Thus, for all small
   enough $\eps$ we have $\kappa(x) \geq 0$ for all $x \in \partial
   \Omega$, which proves the statement.
\end{proof}

The next lemma is key to the analysis of the small $\eps$ regime and
is based on a Bonnesen-type inequality for convex sets with small
isoperimetric deficit \cite{bonnesen24} (for a review, see
\cite{osserman79}). In view of Lemma \ref{lin-isodef}, the latter is
the case for the minimizers of $E_\eps$, when $\eps$ is sufficiently
small. We will use a version of the result that was proved by Fuglede
\cite{fuglede89}, which connects the isoperimetric deficit to the
spherical deviation of the set $\Omega$ from a unit ball centered at
the barycenter of $\Omega$, to prove this lemma.

\begin{lemma}
 \label{lin-fuglede}
  Let $\Omega$ be a minimizer of $E_\eps$, and let $x_0 \in \R^2$ be
  the barycenter of $\Omega$. Then there exists $\eps_3 =
  \eps_3(\alpha) > 0$ such that for all $\eps \leq \eps_3$ 
  \begin{itemize}
  \item[(i)] There exists $\delta > 0$ satisfying    
    \begin{align}
      \label{delD}
      \delta \leq C \sqrt{D(\Omega)},
    \end{align}
   with some universal $C > 0$ such that $B_{1-\delta}(x_0) \subset \Omega
    \subset B_{1+\delta}(x_0)$.
  \item[(ii)] Let $\rho : \R \to (-\delta, \delta)$ be such, that $r =
    1 + \rho(\theta)$ defines the graph of $\partial \Omega$ in polar
    coordinates $(r, \theta)$ centered at $x_0$. Then
    \begin{align}
      \label{fugrho}
      D(\Omega) \leq C || \rho ||^2_{H^1(0, 2 \pi)},
    \end{align}
    for some universal $C > 0$.
  \end{itemize}
\end{lemma}

\begin{proof}
  When $\eps$ is sufficiently small, the minimizer $\Omega$ of $E$
  exists, has small isoperimetric deficit by Lemma \ref{lin-isodef}
  and is convex by Lemma \ref{lin-convex}. The result then follows
  from \cite[Theorem 1.3 and footnote 4]{fuglede89}.
\end{proof}

We can now proceed to the conclusion of the proof of Theorem
\ref{thm-disk}. We start with the far field-dominated case. 

\begin{proposition}[Minimizer is a disk, far field-dominated regime]
  \label{prp-diskfar}
  Let $\alpha < 1$. Then there exists $\eps_0 = \eps_0(\alpha) > 0$,
  such that for all $\eps \leq \eps_0$, the unique, up to
  translations, minimizer of $E_\eps$ is $\Omega = B_1(0)$.
\end{proposition}

\begin{proof}
  If $\eps$ is sufficiently small, there exists a minimizer $\Omega$
  of $E$.  Furthermore, the set $\Omega$ satisfies the conclusions of
  Lemma \ref{lin-fuglede}. Since $\Omega$ is a minimizer, we have
  $E(\Omega) \leq E(B_1(x_0))$, where $x_0$ is the barycenter of
  $\Omega$, which is equivalent to
  \begin{align} \label{dOmballal11} %
    D(\Omega) \leq \frac\eps{2\pi} \left( \int_{B_1(x_0)}
      \int_{B_1(x_0)} {1 \over |x - y|^\alpha} \ dx \, dy -
      \int_\Omega \int_\Omega {1 \over |x - y|^\alpha} \ dx \, dy
    \right).
  \end{align}
  Let $u$ and $u^B$ be the characteristic functions of $\Omega$ and
  $B_1(x_0)$, respectively, and let $v^B$ be as in \eqref{vB}. Then,
  since the non-local kernel is positive-definite (as can be seen from
  \eqref{EOmenonq} and \eqref{Gq}), and since $\int_{\R^2} (u^B - u)
  \, dx = 0$, we have
  \begin{align} 
    \hspace{6ex} & \hspace{-6ex} %
    \int_{B_1(x_0)} \int_{B_1(x_0)} {1 \over |x - y|^\alpha} \ dx \,
    dy - \int_\Omega \int_\Omega {1 \over |x - y|^\alpha} \ dx \, dy
    \notag
    \\
    &= 2 \int_{\R^2} v^B(x - x_0) (u^B(x) - u(x)) \, dx - \int_{\R^2}
    \int_{\R^2} {(u^B(x) - u(x)) (u^B(y) - u(y)) \over |x - y|^\alpha}
    \ dx \, dy
    \notag \\
    & \leq 2 \int_{\R^2} (v^B(x - x_0) - v_0) (u^B(x) - u(x)) \, dx,
 \end{align}
 where $v_0$ is given by \eqref{v0}. Thus
 \begin{align} \label{Enlal1} %
   \hspace{2ex} & \hspace{-2ex} \int_{B_1(x_0)} \int_{B_1(x_0)} {1
     \over |x - y|^\alpha} \ dx \, dy - \int_\Omega \int_\Omega {1
     \over |x - y|^\alpha} \ dx \, dy \notag \\
   & \leq 2 \int_{\Omega \Delta B_1(x_0)} |v^B(x - x_0) - v_0| \, dx
   \leq 2 \ncoL{v^B(\cdot - x_0) - v_0}{\Omega \Delta B_1(x_0)}
   |\Omega \Delta B_1(x_0)| \notag \\
   &\leq C \delta \ncoL{v^B - v_0}{B_{1+\delta}(0) \backslash
     B_{1-\delta}(0)},
  \end{align}
  for some universal $C > 0$.  On the other hand, by Lemma
  \ref{lin-vball}, we have $|v^B - v_0| \leq C \delta$ in
  $B_{1+\delta}(0) \backslash B_{1-\delta}(0)$, for some $C > 0$
  depending only on $\alpha$.  Combining this inequality with by
  \eqref{delD}, \eqref{dOmballal11} and \eqref{Enlal1}, we get
  \begin{align}
    \label{Dal1}
    c \delta^2 \leq D(\Omega) \leq C \eps \delta^2,
  \end{align}
  for some universal $c > 0$ and some $C > 0$ depending only
  $\alpha$. Therefore, as long as $\eps$ is small enough, we have
  $D(\Omega) = 0$, implying that $\Omega = B_1(x_0)$.
\end{proof}

We note that the above proof fails in the near field-dominated regime,
$\alpha \geq 1$, since in this case $v^B$ fails to be in $C^1(\R^2)$,
as can be seen from \eqref{vbrser} (in fact, the radial derivative of
$v^B$ gets singular at $\partial B_1(0)$). Therefore, a more delicate
analysis of the contribution of the deviation of $\Omega$ from a ball
to the non-local part of the energy is necessary. In fact, we need to
prove some cancellations in the difference of the two nonlocal
energies (related to the minimizer and the corresponding ball of the
same area) to obtain an analog of \eqref{Dal1}. For this, we will make
a more detailed use of the Euler-Lagrange equation.

It remains to prove the following proposition.

\begin{proposition}[Minimizer is a disk, near field-dominated regime]
  \label{prp-disknear}
  Let $\alpha \geq 1$. Then there exists $\eps_0 = \eps_0(\alpha) >
  0$ such that for all $\eps \leq \eps_0$, the unique, up to
  translations, minimizer of $E_\eps$ is $\Omega = B_1(0)$.
\end{proposition}

\begin{proof}
  The main point here is to obtain the inequality in the right-hand
  side of (\ref{Dal1}) from (\ref{dOmballal11}). The conclusion then
  follows as in the proof of Proposition \ref{prp-diskfar}. We begin
  by writing
  \begin{align}
    \label{Enlab1}
    \hspace{1ex} & \hspace{-1ex} \int_{B_1(x_0)} \int_{B_1(x_0)} {1
      \over |x - y|^\alpha} \ dx \, dy - \int_\Omega \int_\Omega {1
      \over |x - y|^\alpha} \ dx \, dy
    \notag \\
    & = \int_{B_1(x_0) \backslash \Omega} (v^B(x - x_0) + v(x) - 2
    v_0) \, dx - \int_{\Omega \backslash B_1(x_0)} (v^B(x - x_0) +
    v(x) - 2 v_0) \, dx \notag \\
    & = I + II,
  \end{align}
  where
  \begin{align}
    \label{I}
    I & = \int_{B_1(x_0) \backslash \Omega} (v(x) - v^B(x - x_1(x)))
    \, dx - \int_{\Omega \backslash B_1(x_0)} (v(x) - v^B(x - x_1(x)))
    \, dx, \\
    II & = \int_{B_1(x_0) \backslash \Omega} (v^B(x - x_1(x)) + v^B(x
    - x_0) - 2 v_0) \,
    dx \notag \\
    &- \int_{\Omega \backslash B_1(x_0)} (v^B(x - x_1(x)) + v^B(x -
    x_0) - 2 v_0) \, dx,
    \label{II}
  \end{align}
  and
  \begin{align}
    \label{x1}
    x_1(x) := x_0 + (|x|-1) \, {x - x_0 \over |x - x_0|},
  \end{align}
  i.e., $x_1(x)$ is the center of a ball whose center is shifted from
  $x_0$ in the direction of $x$ in such a way that $x \in \partial
  B_1(x_1(x))$. Introducing polar coordinates as in Lemma
  \ref{lin-fuglede}(ii), we have (with a slight abuse of notation)
  \begin{align}
    \label{Ir}
    I = & \int_0^{2 \pi} \int_1^{1+\rho(\theta)} (
    v^B(r - \rho(\theta)) - v(r, \theta) ) \, r dr d \theta, \\
    \label{IIr}
    II = & \int_0^{2 \pi} \int_1^{1+\rho(\theta)} (2 v_0 - v^B(r) -
    v^B(r - \rho(\theta))) \, r dr d \theta.
  \end{align}

  Let us estimate the term in (\ref{IIr}) first. In the following we
  will only explicitly consider the case $\alpha > 1$, the case
  $\alpha = 1$ is treated analogously. In view of (\ref{vbrser}), with
  $s = r - 1$ we have
  \begin{align}
    \label{IIrrr}
    II & = \int_0^{2 \pi} \int_0^{\rho(\theta)} (2 v_0 - v^B(1+s) -
    v^B(1+s - \rho(\theta))) (1+s) \, ds d \theta \notag \\
    & = C \int_0^{2 \pi} \int_0^{\rho(\theta)} (|s|^{1-\alpha} s -
    |\rho(\theta) - s|^{1-\alpha} (\rho(\theta) - s)) \, ds d \theta +
    O(\delta^2),
  \end{align}
  for some $C > 0$ depending only on $\alpha$, where
 \begin{align}
   \label{del}
   \delta = ||\rho||_{L^\infty(\R)}.    
 \end{align}
 However, the integral in the second line of (\ref{IIrrr}) is
 identically zero, so we have $II = O(\delta^2)$.

  We now turn to estimating (\ref{Ir}), which can be written as
  \begin{align}
    \label{Irrr}
    I = & \int_0^{2 \pi} \int_0^{\rho(\theta)}
    \int_{\theta-\pi}^{\theta + \pi}
    \int_{\rho(\theta')}^{\rho_B(\theta, \theta')} d^{-\alpha}(s, s',
    \theta, \theta') \, (1 + s) (1+s') ds' d \theta' ds d\theta,
  \end{align}
  where $d(s, s', \theta, \theta')$ is the distance between the points
  with polar coordinates $(1 + s, \theta)$ and $(1 + s', \theta')$,
  and $\rho_B(\theta, \theta')$ solves
  \begin{align}
    \label{rhoB}
    1 = (1 + \rho_B)^2 + \rho^2 - 2 \rho (1 + \rho_B) \cos (\theta -
    \theta'), 
  \end{align}
  and for each $\theta$ simply describes the polar graph $r(\theta') =
  1 + \rho_B(\cdot, \theta')$ of a circle shifted by $\rho(\theta)$ in
  the direction of $\theta$ from the origin. Clearly, $\rho_B(\theta,
  \cdot) \in C^\infty(\R)$ for sufficiently small $\delta$, and
  furthermore for all $\theta, \theta' \in \R$ we have
  \begin{align}
    \label{rhoBest}
    |\rho_B(\theta, \theta')| \leq \delta, \qquad |\rho_B(\theta,
    \theta') - \rho(\theta)| \leq C \delta |\theta - \theta'|^2,
  \end{align}
  for some universal $C > 0$. In addition, for small enough $\delta$
  we have
  \begin{align}
    \label{dangles}
  d(s, s', \theta, \theta') \geq c |\theta - \theta'|,  
  \end{align}
  for some universal $c > 0$.

  Combining all the information above, we can write
  \begin{align}
    \label{Ith}
    |I| & \leq C \left| \int_0^{2 \pi} \int_0^{\rho(\theta)}
      \int_{\theta-\pi}^{\theta + \pi}
      \int_{\rho(\theta')}^{\rho_B(\theta, \theta')} |\theta -
      \theta'|^{-\alpha} \, ds' d \theta' ds d\theta \right| \notag \\
    & \leq C \left( \left| \int_0^{2 \pi} \int_0^{\rho(\theta)}
        \int_{\theta-\pi}^{\theta + \pi}
        \int_{\rho(\theta)}^{\rho_B(\theta, \theta')} |\theta -
        \theta'|^{-\alpha} \, ds' d \theta' ds d\theta \right|
    \right. \notag \\
    & \quad ~ \quad + \left. \left| \int_0^{2 \pi}
        \int_0^{\rho(\theta)} \int_{\theta-\pi}^{\theta + \pi}
        \int_{\rho(\theta')}^{\rho(\theta)} |\theta -
        \theta'|^{-\alpha} \, ds' d \theta' ds d\theta \right| \right)
    \notag \\
    & \leq C' \delta \left(\delta + ||\rho_\theta||_{L^\infty(\R)}
      \int_0^{2 \pi} \int_{\theta-\pi}^{\theta + \pi} |\theta -
      \theta'|^{1-\alpha} d \theta' d \theta \right) \notag \\
    & \leq C'' \delta (\delta + ||\rho_\theta||_{L^\infty(\R)}),
  \end{align}
  for some $C, C', C'' > 0$ depending only on $\alpha$, where here and
  below the subscript $\theta$ denotes a derivative with respect to
  $\theta$. So, in order to conclude that $I = O(\delta^2)$ as well,
  it remains to show that
  \begin{align}
    \label{rhoth}
    ||\rho_\theta||_{L^\infty(\R)} \leq C \delta,
  \end{align}
  for some $C > 0$ depending only on $\alpha$.

  To obtain (\ref{rhoth}), we write the Euler-Lagrange equation for
  $\rho(\theta)$ in polar coordinates. Using the well-known formula
  for the curvature in polar coordinates, we can write \eqref{ELeps}
  in the form
  \begin{align}
    \label{ELrho}
    {(1 + \rho)^2 + 2 \rho_\theta^2 - (1 + \rho) \rho_{\theta\theta}
      \over \{ (1 + \rho)^2 + \rho_\theta^2 \}^{3/2} } = {2 \pi \over
      |\partial \Omega|} - 2 \eps (v - \bar v),
  \end{align}
  where $v = v(1 + \rho(\theta), \theta)$. In fact, by continuity of
  $\rho(\theta)$ there exists $\theta^*$ such that $\bar v = v(1 +
  \rho(\theta^*), \theta^*)$. Now, by \cite[Lemma 2.2]{fuglede89}, for
  sufficiently small $\delta$ we also have
  \begin{align}
    \label{rhothinf}
    ||\rho_\theta||_{L^\infty(\R)} \leq C \delta^{1/2},
  \end{align}
  with some
  universal $C > 0$. Therefore, subtracting 1 from both sides of
  (\ref{ELrho}), after a straightforward calculation we obtain
  \begin{align}
    \label{rhoththest}
    ||\rho_{\theta\theta}||_{L^\infty(\R)} \leq C( \delta + D(\Omega)
    + \eps ||v - \bar v||_{L^\infty(\R)}).
  \end{align}
  On the other hand, arguing as in \eqref{Ith}, we have
  \begin{align}
    \label{vrhoth}
    |v(\rho(\theta), \theta) - v_0| & = \left|
      \int_{\theta-\pi}^{\theta+\pi}
      \int_{\rho(\theta')}^{\rho_B(\theta, \theta')}
      d^{-\alpha}(\rho(\theta),
      s', \theta, \theta') (1 + s') \, ds' d \theta' \right| \notag \\
    & \leq C \left( \delta + ||\rho_\theta||_{L^\infty(\R)}
      \int_{\theta-\pi}^{\theta+\pi} |\theta - \theta'|^{1-\alpha} d
      \theta' \right) \notag \\
    & \leq C' (\delta + ||\rho_\theta||_{L^\infty(\R)}),
  \end{align}
  for some $C, C' > 0$ depending only on $\alpha$. In particular,
  since the same estimate holds for $\theta = \theta^*$, we have
  \begin{align}
    \label{vvrhoth}
    ||v - \bar v||_{L^\infty(\R)} \leq C (\delta +
    ||\rho_\theta||_{L^\infty(\R)}), 
  \end{align}
  for some $C > 0$ depending only on $\alpha$.

  Finally, using Lemma \ref{lin-fuglede}(ii), (\ref{rhoththest}) and
  \eqref{vvrhoth}, we conclude that
  \begin{align}
    \label{rhothth2}
    ||\rho_{\theta\theta}||_{L^\infty(\R)} \leq C( \delta + \eps
    ||\rho_\theta||_{L^\infty(\R)}),
 \end{align}
 for some $C > 0$ depending only on $\alpha$. Observe that
 $\rho_\theta(\theta) = \int_{\theta_0}^\theta
 \rho_{\theta\theta}(\theta') d \theta'$ for some $\theta_0 \in
 \R$. Therefore, using smallness of $\eps$, from (\ref{rhothth2}) we
 immediately obtain (\ref{rhoth}).
\end{proof}

\section{Complete characterization in the case of small $\alpha$}
\label{sec-alpha}

In this section, we present the proof Theorem \ref{thm-al0}. The proof
is a slight modification of the proof of Proposition
\ref{prp-diskfar}, and so we find it more convenient to work with the
energy in \eqref{Eeps} (but now without the smallness assumption on
$\eps$). The proof also requires a refinement of the non-existence
result from Sec. \ref{sec:non-exist-minim}.

In terms of $E_\eps$, the result we wish to obtain is a consequence of
the following proposition. 
\begin{proposition} \label{prp-asmall} %
  There exists $\alpha_0 > 0$ such that for all $\alpha \leq \alpha_0$
  the minimizer of $E_\eps$, if it exists, is given by $\Omega =
  B_1(x_0)$, for some $x_0 \in \R^2$.
\end{proposition}

The proof follows from a sequence of lemmas.

\begin{lemma}
  \label{l-Eeps0}
  Let $\Omega \subset \R^2$ be a set of finite perimeter, and let
  $|\Omega| = \pi$. Then
  \begin{align}
    \label{Eeps0}
    E_\eps(\Omega) = |\partial \Omega| + \eps \pi^2 + \alpha
    \int_\Omega \int_\Omega g(x - y) \, dx dy, \qquad |g(x - y)| \leq
    C \eps {\left| \, \ln |x - y| \, \right| \over |x - y|^\alpha},
  \end{align}
  where the constant $C > 0$ depends only on $d := \diam(\Omega)$. 
\end{lemma}

\begin{proof}
  Applying the Taylor formula to the exponential function, we get
  \begin{align}
    \label{g}
   |x - y|^{-\alpha} - 1 = e^{-\alpha \ln |x - y|} - 1 = - \alpha |x - y|^{-
      \alpha \theta} \ln |x - y|,
  \end{align}
  for some $\theta = \theta(x-y) \in (0,1)$. The statement then
  follows with $C = \max \{1, d^2\}$.
\end{proof}

Our next lemma establishes non-existence of minimizers of $E_\eps$ for
sufficiently large $\eps$ {\em uniformly} in $\alpha$ (as long as
$\alpha \leq \alpha_0 < 1$ for some fixed $\alpha_0$).

\begin{lemma}
  \label{l-nonal}
  For every $\alpha_0 \in (0, 1)$ there exists $\eps_2 > 0$ (depending
  only on $\alpha_0$) such that for every $\alpha \in (0, \alpha_0]$
  there is no minimizer of $E_\eps$ for any $\eps > \eps_2$.
\end{lemma}

\begin{proof}
  We prove the statement for the original energy $E$, which amounts to
  existence of $m_2 = m_2(\alpha_0) > 0$ such that there is no
  minimizer of $E$ for all $m > m_2$ and $\alpha \in (0,
  \alpha_0]$. By Lemma \ref{lin-linup}, for a minimizer $\Omega$ of
  $E$ we have $E(\Omega) \leq C m$ for all $m \geq 1$, where the
  dependence of the constant $C > 0$ on $\alpha$ is via
  $E(B_1(0))$. By continuous dependence of $E(B_1(0))$ on $\alpha \in
  [0, \alpha_0]$ (see \eqref{EBR}), we can, in fact, choose $C \geq 1$
  to depend only on $\alpha_0$. Therefore, arguing as in the proof of
  Lemma \ref{lin-d}, we have $m^{2 - \alpha} / C^\alpha \leq m^2 /
  d^\alpha \leq C m$ or, equivalently, $m \leq C^{(1 + \alpha) / (1 -
    \alpha)} \leq C^{2/(1 - \alpha_0)}$.
\end{proof}

We next prove that minimizers of $E_\eps$ must have small isoperimetric
deficit for sufficiently small $\alpha$.

\begin{lemma}
  \label{l-fugal}
  Let $\alpha_0 \in (0, 1)$, let $\alpha \in (0, \alpha_0]$, and let
  $\Omega$ be a minimizer of $E_\eps$. Then $D(\Omega) \leq C \alpha$,
  for some $C > 0$ depending only on $\alpha_0$.
\end{lemma}

\begin{proof}
  By Lemma \ref{l-nonal}, we have $\eps \leq \eps_2(\alpha_0)$, which
  implies, in particular, that $\diam(\Omega) \leq \tfrac12 |\partial
  \Omega| \leq \tfrac12 E(B_1(0)) \leq C$ for some universal $C >
  0$. The result then follows immediately from \eqref{Eeps0} by an
  estimate analogous to the one in \eqref{dEnloc}.
\end{proof}

The result in Lemma \ref{l-fugal} implies that we can use the same
ideas as in Sec. \ref{sec-shape-minim-small} (in the far
field-dominated case), replacing $\eps$ with $\alpha$ and taking
advantage of the smallness of $\alpha$, to prove radial symmetry of
minimizers. In particular, we have the analog of Lemma
\ref{lin-convex}:

\begin{lemma}
  \label{l-convexal}
  There exists $\alpha_0 \in (0, 1)$ such that for every $\alpha \in
  (0, \alpha_0]$ any minimizer $\Omega$ of $E_\eps$ is convex.
\end{lemma}

\noindent Similarly, the analog of Lemma \ref{lin-fuglede} is the
following:
\begin{lemma}
  \label{l-fugal0}
  There exists $\alpha_0 \in (0, 1)$ such that for every $\alpha \in
  (0, \alpha_0]$ any minimizer $\Omega$ of $E_\eps$ satisfies $B_{1 -
    \delta}(x_0) \subset \Omega \subset B_{1+\delta}(x_0)$, where
  $x_0$ is the barycenter of $\Omega$, for some $\delta \leq C
  \sqrt{D(\Omega)}$, with some universal $C > 0$.
\end{lemma}

\begin{proof}[Proof of Proposition \ref{prp-asmall}]
  We argue as in the proof of Proposition \ref{prp-diskfar}. Repeating
  the steps of that proof with the help of Lemmas \ref{l-fugal},
  \ref{l-convexal}, and \ref{l-fugal0}, we obtain
  \begin{align}
    \label{cC}
    c \delta^2 \leq D(\Omega) \leq C \alpha \delta^2,
  \end{align}
 for some universal $c, C > 0$, where the second inequality in \eqref{cC}
  follows from the fact that the potential $v^B$ given by \eqref{vB} obeys
  \begin{align}
    \label{vbral0}
    |\nabla v^B(x)| \leq \alpha \int_{B_1(0)} {1 \over |x -
      y|^{1+\alpha}} \, dy \leq C \alpha \qquad \forall x \in \R^2,
  \end{align}
  for some universal $C > 0$, provided that $\alpha_0$ is sufficiently
  small. The proof is then completed by observing that (\ref{cC})
  implies $D(\Omega) = 0$ for $\alpha_0$ sufficiently small.
\end{proof}

\begin{proof}[Proof of Theorem \ref{thm-al0}]
  Clearly, in view of Proposition \ref{prp-asmall} and Lemma
  \ref{lem-inst-1}(i) there are no minimizers for all $m > m_{c1}$ and
  $\alpha \leq \alpha_0$. It hence remains to show that there exists a
  minimizer for every $m \leq m_{c1}$. The assertion of Theorem
  \ref{thm-al0} then follows by Proposition \ref{prp-asmall}.

  Suppose that $m \leq m_{c1}$ and consider a minimizing sequence $\{
  \Omega_k \}_{k=1}^\infty$. Then $E(\Omega_k) \to e(m) :=
  \inf_{|\Omega| = m} E(\Omega)$ as $k \to \infty$, and by an
  approximation argument we may assume that all sets $\Omega_k$
  consist of $N_k < \infty$ disjoint open connected components. In
  fact, $\Omega_k$ can be chosen so that $N_k$ is independent of
  $k$. Indeed, by Theorem \ref{thm-disk} we can lower the energy by
  replacing all the connected components whose mass is less than $m_0$
  with balls of the same mass translated sufficiently far apart (as in
  the proof of Lemma \ref{prp-connect}). Then, if more than one ball
  is present in the resulting set, by Lemma \ref{lem-inst-1}(ii) we
  can further lower the energy by merging these balls, two at a time,
  and translating the resulting balls further apart.

  In view of the above argument we may assume that $\Omega_k =
  \bigcup_{j=1}^N \Omega_{k}^{(j)}$ with $N_k \leq N$ for some $N \leq
  1 + (m / m_0)$ and the sets $\Omega_k^{(j)}$ are connected (some of
  $\Omega_{k}^{(j)}$ are empty if there are less than $N$ connected
  components). After taking a subsequence, we may assume that for each
  $1 \leq j \leq N$ we have $E(\Omega_k^{(j)}) \to e_j$,
  $|\Omega_k^{(j)}| \to \mu_j$ for some constants $e_j \geq 0$ and
  $\mu_j \geq 0$ as $k \to \infty$, and, furthermore, by compactness
  each set $\Omega_k^{(j)}$ ``converges'' to a set $\Omega^{(j)}$ as
  $k \to \infty$ after a suitable translation. More precisely, if
  $u_k^{(j)}$ are the characteristic functions of $\Omega_k^{(j)}$
  translated to contain the origin, then $u^{(j)}_k \rightharpoonup
  u^{(j)}$ in $BV(\R^2)$ as $k \to \infty$, where $u^{(j)}$ is the
  characteristic function of $\Omega^{(j)}$. Furthermore, since the
  sets $\Omega_k^{(j)}$ are either connected and uniformly bounded or
  empty, we have $|\Omega^{(j)}| = \mu_j$.

  Observe that since $\sum_{j=1}^N \mu_j = m$, we have (see also
  \cite[Remark 4.1]{choksi10})
  \begin{align}
    \label{emless}
    e(m) \leq \sum_{j=1}^N e(\mu_j).
  \end{align}
  Indeed, if $\Omega_j \subset \R^2$ are such that $|\Omega_j| =
  \mu_j$ and $E(\Omega_j) < e(\mu_j) + \delta$ for some $\delta > 0$,
  we can construct a set $\Omega'$ with $|\Omega'| = m$ and
  $E(\Omega') < \sum_{j=1}^N e(\mu_j) + 2 \delta$ by taking $\Omega'$
  to be a union of $\Omega_j$ translated sufficiently far apart. The
  result then follows by arbitrariness of $\delta$. At the same time,
  we have
  \begin{align}
    \label{infEOmj}
    \sum_{j=1}^N E(\Omega^{(j)}) \leq e(m).
  \end{align}
  Indeed, by lower semicontinuity of $E$ with respect to the weak
  $BV$-convergence and by positivity of the kernel in the non-local
  term in the energy we have
  \begin{align}
    \sum_{j=1}^N E(\Omega^{(j)}) \leq \sum_{j=1}^N e_j = \lim_{k \to
      \infty} \sum_{j=1}^N E(\Omega^{(j)}_k) \leq \lim_{k \to \infty}
    E(\Omega_k) = e(m).
  \end{align}
  We now claim that $E(\Omega^{(j)}) = e(\mu_j)$. Indeed, clearly
  $E(\Omega^{(j)}) \geq e(\mu_j)$ for all $1 \leq j \leq N$. On the
  other hand, by \eqref{emless} and \eqref{infEOmj} we get
  \begin{align} \label{dieda} %
    e(m) \ \leq \ \sum_{j=1}^N e(\mu_j) \ \leq \ \sum_{j=1}^N
    E(\Omega^{(j)}) \ \leq \ e(m),
  \end{align}
  so that all inequalities in \eqref{dieda} turn into equalities
  (compare also with \cite[Lemma 4.4(3)]{choksi10}). Again, since
  $e(\mu_j) \leq E(\Omega^{(j)})$ for each $1 \leq j \leq N$, we get
  $e(\mu_j) = E(\Omega^{(j)})$ as well.

  Thus, each set $\Omega^{(j)}$ is a minimizer of $E$ with prescribed
  mass $\mu_j$. Therefore, by Proposition \ref{prp-asmall} for each $1
  \leq j \leq N$ the set $\Omega^{(j)}$ is either a ball or is
  empty. Then, repeating the argument at the beginning of the proof,
  with the help of Lemma \ref{lem-inst-1}(ii) we conclude that
  $E(B_R(0)) \leq \sum_{j=1}^N E(\Omega^{(j)})$, where $R = (m /
  \pi)^{1/2}$, and, hence, $B_R(0)$ is a minimizer by \eqref{infEOmj}.
\end{proof}

\begin{remark}
  It is easy to see that to the leading order in $\alpha$ the
  non-local part of the energy in (\ref{Eeps0}) is generated by the
  kernel $g(x - y) \simeq \eps \ln |x - y|^{-1}$, which appears in the
  studies of the sharp interface version of the Ohta-Kawasaki energy
  in two dimensions \cite{m:cmp10,choksi10,gms11a,topaloglu11}. In
  this respect the result of Proposition \ref{prp-asmall} is closely
  related to the rigidity result obtained in \cite[Proposition
  3.5]{m:cmp10}.
\end{remark}

Finally, let us point out that the fact that the minimizers in Theorem
\ref{thm-al0} exist if and only if $m \leq m_{c1}$, where $m_{c1}$ is
given by \eqref{mc1}, does not rely on smallness of $\alpha$ and would
remain valid as long as minimizers of $E$ are disks whenever they
exist. This can also be seen from the following general result, which
says, essentially that any set of finite perimeter can be replaced by
a set with lower energy consisting of a union of finitely many
disjoint sets, each of which is a minimizer of $E$.

\begin{proposition}
  Let $\Omega \subset \R^2$ be a set of finite perimeter. Then there
  exists a set 
  \begin{align}
    \Omega' = \bigcup_{i = 1}^N \overline\Omega_i, \qquad \quad
    \overline \Omega_i \cap \overline\Omega_j = \varnothing \quad
    \forall i \not= j,
  \end{align}
  with $N < \infty$ and $|\Omega'| = |\Omega|$ such that $E(\Omega')
  \leq E(\Omega)$ and $E(\overline\Omega_i) = \inf_{m = |\Omega_i|}
  E$.
\end{proposition}

\begin{proof}
  If $\Omega$ is a minimizer of $E$, there is nothing to prove. So
  assume that it is not.  Without loss of generality, we may assume
  that $\Omega$ consists of $N < \infty$ disjoint open connected
  components, denoted by $\Omega_i$, and that $\Omega$ has smooth
  boundary. By positivity and decay of the kernel in the non-local
  part of the energy, we have
  \begin{align}
    \label{Einfs}
    E(\Omega) > \sum_{\Omega_i = \overline\Omega_i} \inf_{m = m_i} E +
    \sum_{\Omega_i \not= \overline\Omega_i} \inf_{m = m_i} E, \qquad
    m_i := |\Omega_i|,
  \end{align}
  where $\overline\Omega_i$ are minimizers of $E$ with mass $m_i$,
  whenever such a minimizer exists.  The strict inequality in
  (\ref{Einfs}) follows from the fact that for $N = 1$ the set
  $\Omega$ is not a minimizer, while for $N > 1$ the energy can be
  reduced by spreading different connected components sufficiently far
  apart (as, e.g., in the proof of Lemma \ref{lin-linup}).

  Suppose that $\Omega_i \not= \overline\Omega_i$. If the minimizer of
  $E$ exists for mass $m = m_i$, we replace $\Omega_i$ with $\Omega_i'
  = \overline\Omega_i$. By the minimizing property of
  $\overline\Omega_i$, we then have $E(\overline\Omega_i) = \inf_{m =
    m_i} E$. Alternatively, if the minimum of $E$ is not attained at
  $m = m_i$, there exists $\delta_0 > 0$ such that if $\Omega_i'$ has
  mass $m_i$ and $E(\Omega_i') < \inf_{m = m_i} E + \delta_0$, then
  $\Omega_i'$ is {\em disconnected}.  Indeed, if not, there exists a
  minimizing sequence consisting of $\Omega_k \subset \R^2$ with
  $|\Omega_k| = m_i$ and each $\Omega_k$ connected. Then by the
  argument in the proof of Theorem \ref{thm-exist} the minimum of $E$
  is attained, contradicting our assumption.

  We can, therefore, replace all sets $\Omega_i \not=
  \overline\Omega_i$ for those values of $m_i$ at which the minimum of
  $E$ with $m = m_i$ is not attained with disconnected sets
  $\Omega_i'$ such that
  \begin{align}
    \label{EinfE}
    E(\Omega_i') < \inf_{m = m_i} E + \delta,
  \end{align}
  for some $\delta \in (0, \delta_0)$ to be specified later. Observe
  that by Theorem \ref{thm-exist} for all those we have $m_i > m_1 =
  m_1(\alpha) > 0$. Therefore, the number of the components $\Omega_i
  \not= \overline\Omega_i$ corresponding to $m_i$ at which the minimum
  of $E$ at $m = m_i$ is not attained is bounded above in terms of
  $|\Omega|$. We now apply to each such component the algorithm in the
  proof of Theorem \ref{thm-exist} to lower energy by erasing the
  smallest connected component of $\Omega_i'$ and rescaling the
  resulting set back to mass $m_i$. In view of the fact that $\delta <
  \delta_0$, this process must terminate before only one connected
  component remains. Then, arguing as in the proof of Theorem
  \ref{thm-exist}, we conclude that the mass of each remaining
  connected component of $\Omega_i'$ is bounded below by some $c > 0$
  depending only on $E(\Omega)$ and $\alpha$.

  We are now able to choose $\delta > 0$ sufficiently small and
  construct a new set $\Omega'$ with $|\Omega'| = |\Omega|$ and
  $E(\Omega') < E(\Omega)$ by taking the union of all connected
  components of the sets $\Omega_i'$ constructed above, suitably
  translated to be sufficiently far apart. In this process the mass of
  each connected component of $\Omega'$ that is distinct from a
  minimizer is bounded above by $\max_{\Omega_i \not=
    \overline\Omega_i} |\Omega_i| - c$. Repeatedly applying this
  process, we then find that after finitely many iterations all
  connected components are minimizers.
\end{proof}

\section*{Acknowledgments} 

The authors would like to acknowledge valuable discussions with
R. V. Kohn, M. Novaga and S. Serfaty. C. B. M. was supported, in part,
by NSF via grants DMS-0718027 and DMS-0908279.

\bibliographystyle{plain}

\bibliography{../nonlin,../stat,../mura}

\end{document}